\documentclass{article}
\usepackage[margin=1.75in]{geometry}

\usepackage{amsmath}
\usepackage{amsthm}
\usepackage{amssymb}
\usepackage{bm}
\usepackage{enumitem}
\usepackage{graphicx}
\usepackage{hyperref}
\usepackage{mathtools}
\usepackage{pgf}
\usepackage{stmaryrd}
\usepackage{subfiles}
\usepackage{suffix}
\usepackage{tikz-cd}
\usepackage{titling}
\usepackage[all]{xy}
\xyoption{2cell}
\UseAllTwocells

\theoremstyle{plain}
\newtheorem{thm}{Theorem}[section]

\newtheorem{prop}[thm]{Proposition}
\theoremstyle{definition}
\newtheorem{defn}[thm]{Definition}
\newtheorem{example}[thm]{Example}

\newcommand{\term}[1]{\emph{#1}}

\newcommand{\dash}{\text{--}}
\mathchardef\mhyphen="2D 

\newcommand{\C}{\mathcal{C}}
\newcommand{\D}{\mathcal{D}}
\newcommand{\G}{\mathbb{G}}
\newcommand{\N}{\mathbb{N}}

\newcommand{\Z}{\mathbb{Z}}
\renewcommand{\o}{\circ}
\newcommand{\fo}{\bullet}
\newcommand{\toto}{\rightrightarrows}

\newcommand{\iso}{\cong}
\newcommand{\natto}{\to}
\newcommand{\op}[1]{{#1}^{\mathrm{op}}}

\newcommand{\adjunct}[1]{\widetilde{#1}}
\newcommand{\comma}[2]{(#1 \downarrow #2)}
\renewcommand{\lim}{\operatorname{lim}}
\newcommand{\colim}{\operatorname{colim}}

\newcommand{\catname}[1]{\mathbf{#1}}
\newcommand{\Set}{\catname{Set}}
\newcommand{\Mon}{\catname{Mon}}
\newcommand{\Cat}{\catname{Cat}}
\newcommand{\Alg}[1]{#1\catname{Alg}}
\newcommand{\EndAlg}[1]{#1\catname{Alg}_{\mathrm{end}}}
\newcommand{\MndAlg}[1]{#1\catname{Alg}_{\mathrm{mnd}}}
\newcommand{\End}{\catname{End}}
\newcommand{\Mnd}{\catname{Mnd}}
\newcommand{\FinEnd}{\catname{End}_f}
\newcommand{\FinMnd}{\catname{Mnd}_f}
\newcommand{\Coll}[1]{\catname{Coll}_{#1}}
\newcommand{\Operad}[1]{\catname{Operad}_{#1}}
\newcommand{\CartEnd}{\catname{CartEnd}}
\newcommand{\CartMnd}{\catname{CartMnd}}
\newcommand{\FinCartEnd}{\catname{CartEnd}_f}
\newcommand{\FinCartMnd}{\catname{CartMnd}_f}
\newcommand{\GSet}{\catname{\G Set}}

\newcommand{\Top}{\catname{Top}}
\newcommand{\StrCat}{\catname{Str\mhyphen \omega\mhyphen Cat}}
\newcommand{\ExtN}{\overline{\N}}
\newcommand{\PowerSet}{\mathcal{P}}

\newcommand{\pdpt}{\bullet}

\newcommand{\singleton}{\{ \ast \}}

\DeclareMathOperator{\FreeOperad}{Free}
\DeclareMathOperator{\FreeMnd}{Free}

\newcommand{\del}{\partial}
\newcommand{\yoneda}{\mathbf{y}}
\newcommand{\realize}[1]{|#1|}
\newcommand{\boundary}{i}
\newcommand{\sphere}[1]{S^{#1}}
\newcommand{\disk}[1]{D^{#1}}

\newcommand{\point}{\star}
\newcommand{\id}[1]{\operatorname{id}(#1)}
\DeclarePairedDelimiter{\apply}{[}{]}
\newcommand{\braid}[1]{B_{#1}}

\newcommand\xydiagram[1]{\[ \xymatrix{#1} \]}
\WithSuffix\newcommand\xydiagram*[2]{\[ \xymatrix#1{#2} \]}
\newcommand{\labeledxydiagram}[2]{
  \begin{equation} \label{#1}
  \begin{gathered}
  \xymatrix{#2}
  \end{gathered}
  \end{equation}
}
\newcommand\adjunction[4]{\xymatrix@!{
  {#1} \rtwocell^{#3}_{#4}{'\phantom'\bot} & {#2}
}}
\newcommand\slicearrow[3]{\left( \vcenter{\xymatrix{ #1 \ar[d]_{#2} \\ #3}} \right)}
\newcommand\inlinearrow[3]{\smash{#1 \xrightarrow{#2} #3}}
\newcommand{\pullbackcorner}[1][dr]{\save*!/#1-1.5pc/#1:(-1,1)@^{|-}\restore}
\newcommand{\pushoutcorner}[1][dr]{\save*!/#1+1.5pc/#1:(1,-1)@^{|-}\restore}
\newcommand{\xyhook}{\ar@{^{(}->}}
\newcommand{\xyiso}[1][r]{\ar@{}[#1]|{\iso}}

\begin{document}

\title{Operations in Leinster's Weak $\omega$-Category Operad}
\author{Kyle Raftogianis}
\date{16 May 2017}
\maketitle

\begin{abstract}

In \cite{Batanin98}, Batanin defines a weak $\omega$-category as an algebra for
a certain operad. Leinster refines this idea in \cite{Leinster04}, and defines
the weak $\omega$-category operad as the initial object of a category of
``operads with contraction''. We demonstrate how a higher category structure
arises from this definition by explicitly constructing various composites,
associativity and coherence laws, and an Eckmann-Hilton braiding.

\end{abstract}

\tableofcontents

\section{Introduction}

An $\omega$-category is a generalization of an ordinary category where there are
not only morphisms between objects, but also morphisms between other morphisms.
Its underlying structure consists of a set of \emph{0-cells} (objects),
\emph{1-cells} (morphisms), \emph{2-cells} (transformations), and so on, with
$n$-cells of each dimension $n$. Each $n$-cell has a \emph{source} $(n -
1)$-cell and a \emph{target} $(n - 1)$-cell. We draw low-dimensional cells as
follows (cf.\ \cite[p.\ vi]{Leinster04}):
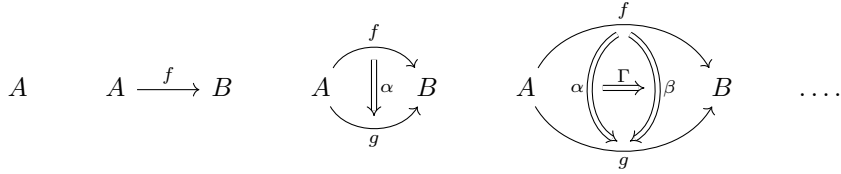
\begin{align*}
  \begin{tikzcd}[ampersand replacement=\&]
    A
  \end{tikzcd}
  &&
  \begin{tikzcd}[ampersand replacement=\&]
    A \ar{r}{f} \& B
  \end{tikzcd}
  &&
  \begin{tikzcd}[ampersand replacement=\&]
    A \arrow[r, bend left=60, "f", ""{below, name=f}]
      \arrow[r, bend right=60, "g"{below}, ""{above, name=g}] \& B
    \arrow[Rightarrow, from=f, to=g, "\alpha"]
  \end{tikzcd}
  &&
  \begin{tikzcd}[ampersand replacement=\&, column sep=0.8in]
    A \arrow[r, bend left=60, "f", ""{below, name=f}]
      \arrow[r, bend right=60, "g"{below}, ""{above, name=g}] \& B
    \arrow[Rightarrow, bend right=60, from=f, to=g, "\alpha"{left}, ""{right, name=alpha}]
    \arrow[Rightarrow, bend left=60, from=f, to=g, "\beta", ""{left, name=beta}]
    \arrow[Rightarrow, from=alpha, to=beta, "\Gamma"]
  \end{tikzcd}
  &&
  \ldots.
\end{align*}
One should imagine the last picture as a 3-dimensional globe, i.e.\ a roughly
spherical shape where $\alpha$ is in front of the page, $\beta$ is behind the
page, and $\Gamma$ goes through the page front-to-back. For this reason, we call
the cellular structure of an $\omega$-category a \emph{globular set}.

An $\omega$-category comes equipped with composition operations that satisfy
associativity and identity laws. If these laws are equalities, the
$\omega$-category is called \emph{strict}; if these laws are only
isomorphisms, the $\omega$-category is called \emph{weak}. For example, in
a strict $\omega$-category we expect to have an associativity equality
\[
  f \o (g \o h) = (f \o g) \o h
\]
whenever this expression makes sense. In contrast, in a weak $\omega$-category
we expect to have only a natural isomorphism
\[
  \alpha_{f,g,h} : f \o (g \o h) \iso (f \o g) \o h
\]
instead of an equality. This natural isomorphism is expected to satisfy
\term{coherence laws} where e.g.\ the \term{pentagon diagram}
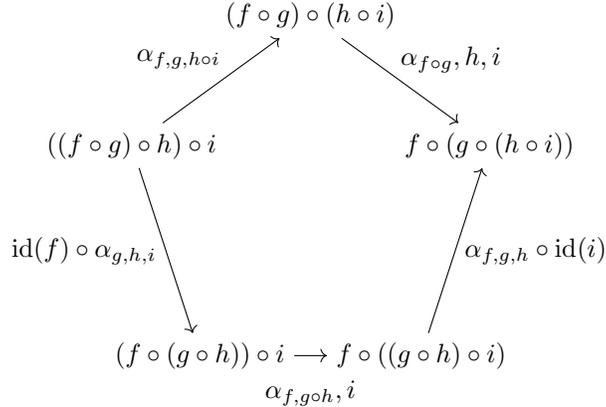
\[
  \begin{tikzpicture}
    \node (P0) at (90:2.5cm) {$(f \o g) \o (h \o i)$};
    \node (P1) at (90+72:2.5cm) {$((f \o g) \o h) \o i$};
    \node (P2) at (90+2*72:2.5cm) {$(f \o (g \o h)) \o i$};
    \node (P3) at (90+3*72:2.5cm) {$f \o ((g \o h) \o i)$};
    \node (P4) at (90+4*72:2.5cm) {$f \o (g \o (h \o i))$};
    \draw
    (P1) edge[->] node[left, near end, outer sep=8pt] {$\alpha_{f,g,h \o i}$} (P0)
    (P0) edge[->] node[right, near start, outer sep=8pt] {$\alpha_{f \o g},h,i$} (P4)
    (P1) edge[->] node[left] {$\id{f} \o \alpha_{g,h,i}$} (P2)
    (P2) edge[->] node[below, outer sep=6pt] {$\alpha_{f, g \o h}, i$} (P3)
    (P3) edge[->] node[right] {$\alpha_{f,g,h} \o \id{i}$} (P4);
  \end{tikzpicture}
\]
commutes. Of course, the two paths in the diagram should not be strictly equal;
rather, we expect to have a natural isomorphism between the two paths that
satisfies its own coherence laws up to natural isomorphism, and so on.

Although strict $\omega$-categories are simpler, weak $\omega$-categories are
much more abundant in nature \cite[p.\ ix]{Leinster04} and therefore more
interesting to study. A motivating example is the \term{fundamental
$\omega$-groupoid} of a topological space $S$, whose 0-cells are points in $S$,
1-cells are paths in $S$, 2-cells are homotopies between paths, and so on. In
addition, weak $\omega$-categories generalize monoidal categories, as evidenced
by the resemblance of the pentagon diagram to the axioms of a monoidal category
\cite[p.\ 162]{MacLane71}. They even generalize \emph{symmetric} monoidal
categories; we show how such symmetry arises from weak structure in Section
\ref{section:eckmann-hilton}. Baez and Dolan have laid out a ``periodic table''
of sorts relating weak categories and various types of monoidal categories
\cite{BaezDolan95}.

Batanin \cite{Batanin98} and Leinster \cite{Leinster04} define these weak
$\omega$-categories in terms of \term{operads}. An operad can be conveniently
characterized as a monad ``lying over'' a base monad, a fact which Leinster
proves as \cite[Corollary 6.2.4]{Leinster04} but otherwise does not play a
significant role in Batanin or Leinster's original work. More recent work on
weak $\omega$-categories uses this characterization of operads \cite{Garner10,
Batanin11}. We take this relationship between operads and monads further by
showing that free operads are in fact free monads, and this monadic
characterization plays a central role in our construction of various weak
$\omega$-category operations.

A \term{globular operad} is a monad lying over the monad for strict
$\omega$-categories. Intuitively, a globular operad consists of a collection of
``operations'', and a system for composing these operations. An \term{algebra}
for a globular operad consists of a system of ``applying'' these operations. For
example, the operad for strict $\omega$-categories has operations that are
``pasting diagrams'', or strict formal composites of $\omega$-category cells
(not to be confused with composites of operations). A strict operation in a
strict $\omega$-category might look like
\[
  \begin{tikzcd}
    A \arrow[r, bend left=75, looseness=2, "f", ""{below, name=f}]
      \arrow[r, bend left=35, "g"{above, near start}, ""{above, name=g1}, ""{below, name=g2}]
      \arrow[r, bend right=35, "h"{below, near start}, ""{above, name=h1}, ""{below, name=h2}]
      \arrow[r, bend right=75, looseness=2, "i"{below}, ""{above, name=i}]
    & B
      \arrow[r, "j"]
    & C
      \arrow[r, bend left, "k", ""{below, name=k}]
      \arrow[r, bend right, "l"{below}, ""{above, name=ell}]
    & D
    \arrow[Rightarrow, from=f, to=g1, "\alpha"]
    \arrow[Rightarrow, from=g2, to=h1, "\beta"]
    \arrow[Rightarrow, from=h2, to=i, "\gamma"]
    \arrow[Rightarrow, from=k, to=ell, "\delta"]
  \end{tikzcd}
\]
where $\alpha, \beta, \ldots$ are 2-cells, $f, g, \ldots$ are 1-cells, and $A,
B, \ldots$ are 0-cells, and the 2-cells are pasted together along their
boundaries. An algebra for this operad is a strict $\omega$-category, so in a
strict $\omega$-category we may ``apply'' this 2-dimensional operation to obtain
a composite 2-cell.

Leinster defines a weak globular operad to be one equipped with a
\term{contraction}, which specifies a way to ``lift'' strict operations into
weak operations. The associativity and identity isomorphisms of a weak globular
operad come from lifting strict associativity and identity operations. We think
of the operations of a weak globular operad as ``weak composites'' of cells, so
that an algebra for such an operad is a system of applying these weak operations
to form composite cells. A weak $\omega$-category is by definition an algebra
for a weak globular operad.

This paper is organized as follows. Section \ref{section:monads} begins with a
brief introduction to monads and some building blocks that operad theory rests
upon, namely monad morphisms, finitary monads, and cartesian monads. We review
how to construct a free monad on an endofunctor, and we prove that the free
monad of a finitary and cartesian endofunctor is also finitary and cartesian. In
Section \ref{section:operads} we introduce the generalized operads that are the
subject of Leinster's book \cite{Leinster04}, with an emphasis on their relation
to monads. We discuss the prototypical example of \emph{planar operads}, which
are operads consisting $n$-ary operations for $n \in \N$; additionally, we show
that the free operad construction is simply an instance of the free monad
construction. We describe globular operads in Section
\ref{section:globular-operads}, and review Leinster's definition of a
contraction. Finally, in Section \ref{section:weak-operations} we demonstrate
how weak $\omega$-category structure arises from such a contraction. We first
use the contraction structure to compose cells in $\omega$-category, which takes
some care. This extra effort is worth it, however, as we show that using the
same contraction structure to define both composites and associativity laws
makes constructing the required coherence and naturality laws easy. Finally, we
invoke the Eckmann-Hilton argument to construct a ``braiding'' that rotates two
cells around each other, demonstrating how symmetry can spontaneously arise from
weak structures.

\subsection*{Acknowledgements}

This paper was written as an undergraduate thesis at the University of
California, Berkeley, advised by Professor Thomas Scanlon. I am grateful for his
time, valuable advice, and encouragement to explore.

\section{Monads} \label{section:monads}

A \term{monad} on a category $\C$ is a structure that encodes an ``algebraic
theory''. Associated to each monad $T$ is a category $\Alg{T}$ of
``algebras'' of that theory, each of which consists of an object
of $\C$ endowed with ``structure'' (i.e.\ operations satisfying laws) according
to the theory that $T$ represents. For example, there is a monad
$T_\catname{Grp}$ on $\Set$ that represents the theory of groups, so that
$\Alg{T_\catname{Grp}}$ is the category of groups. Indeed, any category of
algebraic structures defined by an \term{equational presentation}
induces a monad whose algebras are the structures presented,
as we sketch in Example \ref{example:equational-presentation}.

We can define richer algebraic theories by constructing monads on categories
other than $\Set$. For example, a category can be thought of as a graph with
composition and identity operations, so a category is an algebra for some monad
on the category of graphs. We define weak $\omega$-categories as algebras
for a certain monad.

\subsection{Monads and Algebras}

We begin by reviewing some basic definitions, and giving a few examples of
monads. One can refer to Awodey \cite{Awodey11} or Mac Lane \cite{MacLane71} for
a more thorough account of the subject of monads.

\begin{defn}
  A \term{monad} $(T, \eta, \mu)$ on a category $\C$ consists of an
  endofunctor $T : \C \to \C$, a natural transformation $\eta : 1_\C \natto T$
  (the \term{unit}), and a natural transformation $\mu : T^2 \natto T$ (the
  \term{multiplication}) such that the following diagrams of natural
  transformations commute.
  \begin{align*}
    \xymatrix{
      T^3 \ar[r]^{T \mu} \ar[d]_{\mu T} & T^2 \ar[d]^{\mu} \\
      T^2 \ar[r]_{\mu} & T
    }
    & &
    \xymatrix{
      T \ar[rd]_{1} \ar[r]^{\eta T} & T^2 \ar[d]^-{\mu}
        & T \ar[l]_{T \eta} \ar[ld]^{1} \\
      & T
    }
  \end{align*}
  We will often refer to a monad $(T, \eta, \mu)$ as simply $T$.
\end{defn}

Every pair of adjoint functors $F \dashv U$ induces a monad on $\C$ where
\begin{itemize}
  \item the underlying endofunctor is $T = UF : \C \to \C$,
  \item the unit $\eta$ is the unit of the adjunction $\eta : 1 \natto UF$, and
  \item the multiplication $\mu$ is the natural transformation
  \[ \mu = U\epsilon F : UFUF \natto UF \]
  using the counit $\epsilon : FU \natto 1$ of the adjunction.
\end{itemize}
We draw such an adjunction $F \dashv U$ as
\[ \adjunction{\C}{\D}{F}{U} \]
where the tip of the $\dashv$ points to the left adjoint. The induced
monad $T$ is the clockwise composite $\C \to \C$.

In fact, every monad on $\C$ is induced in this way from a canonical adjunction,
where $\D$ is the following category.

\begin{defn}
  Let $T$ be a monad on the category $\C$. The category $\Alg{T}$ of
  \term{$T$-algebras} is the category where
  \begin{itemize}
    \item objects are pairs $(A, f)$ consisting of an object $A \in \C$ and
    a map $\inlinearrow{TA}{f}{A}$
    that commutes with the unit and multiplication, meaning that the following
    diagrams commute.
    \begin{align*}
      \xymatrix{
        A \ar[rd]_{1} \ar[r]^{\eta_A} & TA \ar[d]^{f} \\
        & A
      }
      & &
      \xymatrix{
        T^2A \ar[d]_{\mu_A} \ar[r]^{Tf} & TA \ar[d]^{f} \\
        TA \ar[r]_{f} & A
      }
    \end{align*}
    \item a morphism
    $\left( \inlinearrow{TA}{f}{A} \right) \to \left( \inlinearrow{TB}{g}{B} \right)$
    is a map $h : A \to B$ such that the diagram
    \xydiagram{
      TA \ar[d]_f \ar[r]^{Th} & TB \ar[d]^g \\
      A \ar[r]_h & B
    }
    commutes.
  \end{itemize}
\end{defn}

We call an object of $\Alg{T}$ a \term{$T$-algebra} and a morphism of
$\Alg{T}$ a \term{$T$-algebra homomorphism}.
The category $\Alg{T}$ is also called the \term{Eilenberg-Moore category} of $T$
after the following theorem, which establishes the canonical adjunction that
induces $T$.

\begin{thm}[Eilenberg-Moore]
  Let $T$ be a monad on $\C$, and let
  $U : \Alg{T} \to \C$ be the forgetful functor that sends an algebra
  $\inlinearrow{TA}{f}{A}$ to its underlying object $A$.
  Then $U$ has a left adjoint
  \[ \adjunction{\C}{\Alg{T}}{F}{U} \]
  and the monad induced by this adjunction is the monad $T$.
\end{thm}
\begin{proof}
  The left adjoint $F$ sends an object $A$ to the \term{free algebra}
  $\inlinearrow{T^2A}{\mu_A}{TA}$.
  We take the unit $\eta$ of the adjunction to be the unit of monad,
  and the counit to be the natural transformation whose
  component at an algebra $\inlinearrow{TA}{f}{A}$ is the algebra morphism
  \xydiagram{
    T^2A \ar[d]_{Tf} \ar[r]^{\mu_A} & TA \ar[d]^{f} \\
    TA \ar[r]_{f} & A
  }
  which commutes by the definition of an algebra.
\end{proof}

We say a functor $U$ is \term{monadic} if it is canonically equivalent
to the forgetful functor from a category of algebras. Let
\[ \adjunction{\C}{\D}{F}{U} \]
be an adjunction, and let $T = UF$ be its induced monad.
The canonical \term{comparison functor} $K^T : D \to \Alg{T}$ is given
on objects by
\[
  K^T(A) = \slicearrow{TUA}{U\epsilon_A}{UA}.
\]

\begin{defn}
  A functor $U : \D \to \C$ is \term{monadic} if $U$ has a left adjoint and
  the induced comparison functor $K^T$ is an equivalence.
\end{defn}

In other words, if $U : \D \to \C$ is monadic, then it induces a monad $T$ on
$\C$, and we may think of $\Alg{T}$ as just being the category $\D$.

\begin{example} \label{example:free-monoid-monad}
  A \term{monoid} is an algebraic structure consisting of a set $X$, an
  element $e \in X$, and a binary operation $\cdot : X \times X \to X$ so that
  for all $x, y, z \in X$
  \[
    \begin{array}{ccc}
      e \cdot x = x & x \cdot e = x & x \cdot (y \cdot z) = (x \cdot y) \cdot z
    \end{array}
  \]
  i.e.\ a monoid is ``group without inverses''. Let $\Mon$ be the category of
  monoids and monoid homomorphisms.

  The forgetful functor $U : \Mon \to \Set$ that projects
  out the underlying set has a left adjoint
  \[ \adjunction{\Set}{\Mon}{F}{U}. \]
  The underlying set of $FX$ (that is, the set $TX = UFX$) is
  \[
    TX = \coprod_{n \in \N} X^n
  \]
  so that an element of $TX$ is an $n$-tuple $(x_1, \ldots, x_n)$ for some
  $n$, a \term{list} of length $n$. The monoid operation
  $\cdot : TX \times TX \to TX$ is the concatenation
  \[
    (x_1, \ldots, x_m) \cdot (y_1, \ldots, y_n)
      \mapsto (x_1, \ldots, x_m, y_1, \ldots, y_n)
  \]
  and the identity element is the empty list $() \in TX$. The unit
  $X \to TX$ injects an $x \in X$ as the one-element list $(x)$.

  We may think of a list as a ``formal composite'' of elements of $X$. For a
  monoid $M$, the counit $FUM \to M$ is the monoid homomorphism that evaluates
  this formal composite as
  \begin{align*}
    () &\mapsto e \\
    (m_1, m_2, \ldots, m_n) &\mapsto m_1 \cdot m_2 \cdot \ldots \cdot m_n.
  \end{align*}
  It follows that the multiplication $\mu : T^2 X \to TX$ of the monad $T$
  concatenates a nested list, e.g.
  \[
    ((x_1, x_2, x_3), (x_4), (x_5, x_6)) \mapsto (x_1, x_2, x_3, x_4, x_5, x_6).
  \]

  The comparison functor $K^T : \Mon \to \Alg{T}$ uses the counit to turn a
  monoid $M$ into a $T$-algebra $TM \to M$. Conversely, there is a functor
  $\Alg{T} \to \Mon$ which takes an algebra $\inlinearrow{TX}{f}{X}$ to the
  monoid $(X, e, \cdot)$ where $e$ is the result of applying $f$ to the empty
  list, and $x \cdot y$ is the result of applying $f$ to the two-element list
  $(x, y)$. This functor is inverse to the comparison functor, so the
  forgetful functor $U : \Mon \to \Set$ is monadic, and thus algebras for the
  monad $T$ are equivalent to monoids.
\end{example}

\begin{example} \label{example:equational-presentation}
  Categories of algebraic structures are typically defined by an
  \term{equational presentation}, which says the structures being defined are
  ``sets with operations satisfying axioms''. These include
  the categories of monoids, groups, and rings, although notably not the
  category of fields. Here we only consider \term{finitary} equational
  presentations, meaning we have finitely many operations and finitely many
  laws.

  A \term{finitary signature} is a set
  $\Sigma$ together with a function $\Sigma \to \N$. The elements of
  $\Sigma$ are called \term{function symbols} and the function $\Sigma \to \N$
  assigns each function symbol to its \term{arity}.
  An \term{algebra} for $\Sigma$ consists of a set $X$ equipped with, for
  each function symbol $\sigma \in \Sigma$ of arity $n$, a function
  $\sigma_X : X^n \to X$.
  A \term{homomorphism} of algebras $X$ and $Y$ is a function $X \to Y$
  such that for each operation $\sigma$ of arity $n$, the diagram
  \xydiagram{
    X^n \ar[d]_{\sigma_X} \ar[r]^{f^n} & Y^n \ar[d]^{\sigma_Y} \\
    X \ar[r]_{f} & Y
  }
  commutes. Thus the algebras of a signature form a category, which we
  denote $\Alg{\Sigma}$.

  For example, to define an abelian group in this way, we start with three
  function symbols $\Sigma = \{ +, -, 0 \}$ of arities 2, 1, and 0 respectively.
  A $\Sigma$-algebra is a set $X$ equipped with a binary operation
  $+_X : X \times X \to X$, a unary operation $-_X : X \to X$,
  and a nullary operation $0_X : 1 \to X$ (i.e.\ an element $0_X \in X$).
  A homomorphism between $\Sigma$-algebras $X$ and $Y$ is a function $f : X \to Y$
  such that for all $x, y \in X$,
  \begin{align*}
    f(x + y) &= f(x) + f(y) \\
    f(-x) &= -f(x) \\
    f(0_X) &= 0_Y.
  \end{align*}

  Of course, the definition of an abelian group would not be complete without
  axioms. We also form a finite set $E$ of \term{equations}, which are pairs of formal
  expressions involving formal variables and operations $\sigma \in \Sigma$
  (more rigorously, pairs of elements of the free $\Sigma$-algebra on some
  fixed countable set of variables $V = \{x, y, z, \ldots \}$). For abelian
  groups, this set of equations would be
  \begin{align*}
    (x + y) + z &= x + (y + z) & x + y &= y + x \\
    (-x) + x &= 0 & 0 + x &= x \\
    x + (-x) &= 0 & x + 0 &= x.
  \end{align*}
  A \term{finitary equational presentation} is a pair $(\Sigma, E)$ of a finitary signature and
  finitely many equations, and the category of algebras $\Alg{(\Sigma, E)}$ for such an equational
  presentation is the full subcategory of $\Alg{\Sigma}$ comprising $\Sigma$-algebras where the equations are
  satisfied. For the example above, the category of
  algebras $\Alg{(\Sigma, E)}$ is precisely the category of abelian groups.

  For every equational presentation $(\Sigma, E)$, we have a forgetful functor
  $\Alg{(\Sigma, E)} \to \Set$ that sends a $\Sigma$-algebra to its underlying
  set. This functor is monadic \cite{Adamek94}, meaning that we have a left
  adjoint $\Set \to \Alg{(\Sigma, E)}$ and the category of algebras for the
  induced monad $T_{(\Sigma, E)} : \Set \to \Set$ is equivalent to
  $\Alg{(\Sigma, E)}$. In our example above, there is a monad
  $T_{\catname{Ab}} : \Set \to \Set$ whose algebras are abelian groups, and
  this monad takes a set $X$ to the set underlying the free abelian group on
  $X$, i.e.\ finite $\Z$-linear combinations of elements of $X$.

  This result justifies our thinking of monads as
  ``generalized algebraic theories'', and is the origin of the terms
  ``algebra'' and ``homomorphism'' for the objects and morphisms of an
  Eilenberg-Moore category.
\end{example}

\begin{example} \label{example:free-category-monad}

As a final example,
we construct a monad whose algebras are ordinary 1-categories. Let
$\catname{Grph}$ be the category of (directed multi-)graphs, where an object of
$\catname{Grph}$ consists of a set $E$ of edges and a set $V$ of vertices, with
two functions $E \toto V$ that assigns each for edge
to its start and end vertices respectively. More precisely, the category of
$\catname{Grph}$ is the presheaf category $[\op{Q}, \Set]$, where
$Q$ is the category $0 \toto 1$. Let $T$ be the functor which on
objects takes a graph $G$ to a graph with the same vertices $G$,
but whose edges are \emph{paths} of consecutive edges in $G$, as in the
following picture.
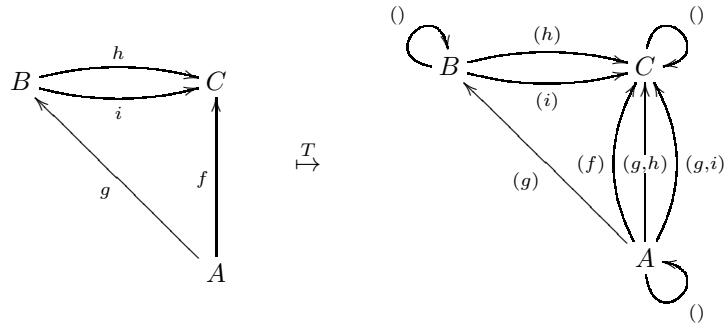
\[
  \vcenter{\xymatrix@=60pt{
    B \ar@/^/[r]^h \ar@/_/[r]_i & C \\
    & A \ar[u]^f \ar[ul]^g
  }} \qquad \overset{T}{\mapsto} \qquad
  \vcenter{\xymatrix@=60pt{
    B \ar@(l,u)^{()} \ar@/^/[r]^{(h)} \ar@/_/[r]_{(i)} & C \ar@(u,r)^{()} \\
    & A \ar@(d,r)_{()} \ar[ul]^{(g)} \ar@/^1pc/[u]^{(f)} \ar[u]|{(g, h)} \ar@/_1pc/[u]_{(g, i)}
  }}
\]

An algebra for this monad will consist of a graph $C$ (whose vertices and edges
we think of as objects and morphisms respectively) together with a graph
homomorphism $TC \to C$ that takes a path of morphisms to a composite morphism
in $C$. This operation must respects the unitality and associativity of
composition of paths, so an algebra for $T$ is precisely a category.

Note that an $T$-algebra homomorphism is exactly a functor of categories. An
isomorphism in this category is a strict isomorphism consisting of functors $F$
and $G$ such that $FG = 1$ and $GF = 1$. However, this is not the correct notion
of equivalence for categories, which requires instead that there only be natural
isomorphisms $FG \iso 1$ and $GF \iso 1$. Similarly, defining a weak
$\omega$-category as an algebra for some monad is not enough to discuss when two
weak $\omega$-categories are only weakly equivalent.

\end{example}

\subsection{Categories of Monads}

For a category $\C$,
endofunctors $\C \to \C$ form a category $\End(\C)$ whose morphisms are
natural transformations of endofunctors. Similarly, we can define a category
$\Mnd(\C)$ of monads on $\C$ whose morphisms are \term{monad morphisms},
which we now define.
Street \cite{Street72} gives a more general notion of a monad morphism than
defined here that relates monads on different categories, but we do not need its
full generality.

\begin{defn}
If $S$ and $N$ are both monads on $\C$, a \term{monad morphism} $S \to T$ is a
natural transformation $\alpha : S \natto T$ that commutes with the unit and
multiplication of the monads $S$ and $T$, meaning that the following diagrams
commute.
\begin{align*}
  \xymatrix{
    & 1 \ar[ld]_{\eta} \ar[rd]^{\eta} \\
    S \ar[rr]_{\alpha} & & T \\
  } & &
  \xymatrix{
    S^2 \ar[r]^{\alpha \alpha} \ar[d]_{\mu} & T^2 \ar[d]^{\mu} \\
    S \ar[r]^{\alpha} & T \\
  }
\end{align*}
\end{defn}

Monad morphisms are closed under composition, so we may form a category
$\Mnd(\C)$ of monads on $\C$ and monad morphisms, and we have a faithful
forgetful functor $\Mnd(\C) \to \End(\C)$.

Monad morphisms are characterized by their action on categories of algebras.
Every monad morphism $\alpha : S \natto T$ induces a functor
\[
  \begin{array}{cccc}
    \alpha^* : &\Alg{T} &\to &\Alg{S} \\
      &\slicearrow{TA}{f}{A} &\mapsto
      & \left( \vcenter{\xymatrix{ SA \ar[d]_{\alpha_A} \\ TA \ar[d]_{f} \\ A}} \right)
  \end{array}
\]
between categories of algebras in the opposite direction, lying over $\C$ in the
sense that the diagram of functors
\[
\xymatrix@=15pt{
  \Alg{T} \ar[ddr]_U \ar[rr]^{\alpha^*} & & \Alg{S} \ar[ddl]^U \\
  \\
  & \C
}
\]
commutes. In fact, we can show that any functor of algebras lying over $\C$ is
$\alpha^*$ for some monad morphism $\alpha$.

\begin{prop} \label{prop:star-inverse}
  Let $S$ and $T$ be monads on the category $\C$. Then the function $(\dash)^*$
  is a bijection
  \[
    (\dash)^* : \Mnd(\C)(S, T) \iso \Cat/\C(\Alg{T}, \Alg{S}).
  \]
\end{prop}

The proof is specialized from a more general proof that Leinster gives
as \cite[Lemma 6.1.1]{Leinster04}.

\begin{proof}
  Suppose we are given a functor $F : \Alg{T} \to \Alg{S}$ lying over $\C$.
  For any object $A \in \C$, the functor $F$ sends the free $T$-algebra
  $\inlinearrow{T^2A}{\mu}{TA}$
  to some $S$-algebra
  $\inlinearrow{STA}{\phi_A}{TA}$
  lying over the same object $TA$. We define a component $\alpha_a : SA \to TA$
  by the composite
  \xydiagram{
    SA \ar[r]^{S\eta_A} & STA \ar[r]^{\phi_A} & TA.
  }
  These components assemble into a natural transformations
  $\alpha_F : S \natto T$, and we may verify using the algebra laws that this
  natural transformation is a monad morphism. The function $F \mapsto \alpha_F$
  is the desired inverse to $(\dash)^*$.
\end{proof}

Thus in some sense giving a monad morphism $S \to T$ is the same
making the statement ``every $T$-algebra is canonically an $S$-algebra''.

\begin{example}
  Let $T_\Mon$ be the free monoid monad from Example
  \ref{example:free-monoid-monad}, and let $T_{\catname{Grp}}$ be the monad
  whose algebras are groups (which can be defined by Example \ref{example:equational-presentation}).
  Any group is also a monoid by forgetting
  inverses, so we have a functor $\catname{Grp} \to \Mon$ (that is,
  $\Alg{T_{\catname{Grp}}} \to \Alg{T_\Mon}$) lying over the forgetful functors
  to $\Set$. Thus we have a monad morphism $T_\Mon \to T_{\catname{Grp}}$, whose
  component at $X$ interprets a list in the free monoid on $X$ as a word in the
  free group on $X$.

  Note however that the functor $(\dash)^\times : \catname{Rng} \to \catname{Grp}$
  that sends a ring $R$ to its group of units $R^\times$ is \emph{not} induced
  by any monad morphism, as if $R$ is nontrivial then the underlying set of
  $R^\times$ is a proper subset of $R$.
\end{example}

\subsection{Finitary Monads} \label{section:finitary-monads}

A monad on $\Set$ is induced by a finitary equational
presentation if and only if it is a \term{finitary monad} \cite{Adamek94},
in a sense we now define. We first need the following technical definition.

\begin{defn}
  A small category $J$ is \term{filtered} if every finite diagram in $J$
  has a cocone (or ``upper bound'') in $J$.

  A \term{filtered colimit} in a category $\C$ is a colimit over a diagram
  $J \to \C$ where $J$ is a filtered category.
\end{defn}

An example of a filtered category that we use extensively is the category
\xydiagram{
  0 \ar[r] & 1 \ar[r] & 2 \ar[r] & \ldots.
}
A colimit in $\Set$ over this category is a union of increasing subsets.

Limits of type $P$ and colimits $J$ are said to \emph{commute}
in a category $\C$ if for any diagram $F : P \times J \to \C$ the canonical map
\[
  \colim_J \lim_P F \to \lim_P \colim_J F
\]
is an isomorphism.
A crucial property of filtered colimits is that they commute with finite
limits in $\Set$ \cite[\S IX.2]{MacLane71}, and hence in any presheaf
category.
It follows that any functor on a presheaf category defined by colimits and
finite limits will commute with (i.e.\ preserve) filtered colimits. We call
such a functor \term{finitary}.

\begin{defn}
  Let $\C$ be a presheaf category.
  A functor $\C \to \C$ is finitary if it preserves filtered colimits.
  A monad $T$ is finitary if its underlying endofunctor $T : \C \to \C$ is
  finitary.
\end{defn}

For example, the free monoid monad
\[
  TX = \coprod_{n \in \N} X^n
\]
is finitary, as it is the coproduct of finite limits $X^n$.

\subsection{Cartesian Monads}

\begin{defn}
  We call a structure having the property of ``preserving pullbacks''
  \term{cartesian}.

  \begin{itemize}
  \item A category is cartesian if it has all pullbacks.%
  \footnote{This differs from the usual usage of the term ``cartesian'', which
  usually means that the category has finite (cartesian) products. A category having all
  pullbacks is more accurately called \emph{locally} cartesian, since a pullback is
  a product in a slice category.}

  \item A functor is cartesian if it preserves pullbacks i.e.\ the image of
  a pullback square is another pullback square.

  \item A natural transformation
  $\alpha : F \natto G$ is
  cartesian if for all $f : X \to Y$ the naturality square
  \xydiagram{
    FX \ar[r]^{Ff} \ar[d]_{\alpha_Y} & FY \ar[d]^{\alpha_Y} \\
    GX \ar[r]_{Gf} & GY
  }
  is a pullback square.

  \item A monad $(T, \eta, \mu)$ on a cartesian category $\C$ is cartesian if
  the functor $T$ is cartesian and
  and the natural transformations $\eta$ and $\mu$ are cartesian.

  \item A monad morphism $\alpha : S \natto T$ is cartesian if its underlying
  natural transformation is cartesian.
  \end{itemize}
\end{defn}

Pullbacks commute with coproducts in $\Set$, and thus in any presheaf category.
It follows that any functor defined by coproducts and
limits on a presheaf category will be cartesian.

For any category $\C$, we may form a category $\CartEnd(\C)$ of
cartesian endofunctors $\C \to \C$ and cartesian natural transformations. If
$\C$ is a presheaf category,
pullbacks in $\C$ commute with limits, coproducts, and filtered colimits. In
this case, the faithful functor $\CartEnd(\C) \to \End(\C)$ creates limits,
coproducts, and filtered colimits. For example, if $F, G, H : \C \to \C$ are
functors, the coproduct inclusions $F \to F + G$ and $G \to F + G$
are cartesian natural transformations, and if $\alpha : F \to H$ and
$\beta : G \to H$ are cartesian then so is the natural transformation
$[\alpha, \beta] : F + G \to H$.

Cartesianness is a strong property to impose on a natural transformation.
If the codomain of $F$ and $G$ has a terminal object, then a cartesian natural
transformation $\alpha : F \to G$ is determined by the component
$\alpha_1 : F1 \to G1$, as any other component $\alpha_A$ is the pullback
along $\inlinearrow{GA}{G!}{G1}$.
\xydiagram{
  FA \ar[r]^{F!} \ar[d]_{\alpha_A} \pullbackcorner & F1 \ar[d]^{\alpha_1} \\
  GA \ar[r]_{G!} & G1
}
Furthermore, many properties of $G$ can be ``pulled back'' by $\alpha$ to $F$.

\begin{prop} \label{prop:cartesian-list}
  Let $\C$ and $\D$ be presheaf categories, let $F, G : \C \to \D$ be functors,
  and let $\alpha : F \to G$ be a cartesian natural transformation.

  \begin{enumerate}
    \item If $G$ is cartesian, then $F$ is cartesian also.
    \item If $G$ is finitary, then $F$ is finitary also.\footnote{%
    Garner mentions this fact in passing in \cite{Garner10}.}
  \end{enumerate}
\end{prop}
\begin{proof} \leavevmode
  \begin{enumerate}
    \item Suppose
    \xydiagram{
      A \ar[r] \ar[d] \pullbackcorner & B \ar[d] \\
      C \ar[r] & D
    }
    is a pullback square in $\C$. Applying $\alpha$ we may form the commuting cube
    \[
    \xymatrix@!0{
      & FA \ar[rr]\ar'[d][dd] \ar[dl]
      & & FB \ar[dd] \ar[dl]
      \\
      FC \ar[rr]\ar[dd]
      & & FD \ar[dd]
      \\
      & GA \ar'[r][rr] \ar[dl]
      & & GB \ar[dl]
      \\
      GC \ar[rr]
      & & GD
      }
    \]
    where the lower 5 faces are pullback squares. Thus the top
    face is a pullback also.

    \item Since $\D$ has a terminal object, for any $A$ we have a pullback
    square
    \xydiagram{
      FA \ar[r]^{F!} \ar[d]_{\alpha_A} \pullbackcorner & F1 \ar[d]^{\alpha_1} \\
      GA \ar[r]_{G!} & G1.
    }
    Suppose $D : J \to \C$ is a diagram where $J$ is a filtered category.
    Let $P$ be the category whose limits are pullbacks, namely
    \xydiagram{
      & \cdot \ar[d] \\
      \cdot \ar[r] & \cdot
    }
    and define the composite diagram $D' : P \times J \to \C$ where the component
    $P \to \C$ at $j \in J$ is
    \xydiagram{
      & F1 \ar[d]^{\alpha_1} \\
      G(j) \ar[r]_{G!} & G1.
    }
    Then by commutativity of filtered colimits and finite limits
    \[
      F (\colim_J D) \iso \lim_P \colim_J D' \iso \colim_J \lim_P D' \iso \colim_J FD.
    \]
  \end{enumerate}\vspace*{-2\baselineskip}
\end{proof}

\begin{example} \label{free-monoid-monad-cartesian}
  The free monoid monad is cartesian. Its functor part
  \[
    TX = \coprod_{n \in \N} X^n
  \]
  is a coproduct of limits, and so is cartesian.

  Leinster shows that the unit and multiplication of $T$ are cartesian by
  direct computation in \cite{Leinster98}; here we sketch a more abstract
  approach. The unit is the coproduct
  inclusion $X \to \coprod_{n \in \N} X^n$ at $n = 1$, and thus is cartesian.
  Again using the commutativity of pullbacks and coproducts, showing that
  $\mu : T^2 X \to TX$ is cartesian reduces to showing the transformations
  \begin{align*}
    e_X : 1 \to TX & & (\dash \cdot_X \dash) : TX \times TX \to TX
  \end{align*}
  are cartesian. This in turn reduces to showing that
  \begin{align*}
    1 \to X^0 & & X^m \times X^n \to X^{m + n}
  \end{align*}
  are cartesian. But these are natural isomorphisms and trivially
  cartesian, so $\mu$ is also.
\end{example}

Lastly, we characterize which functors $\Alg{T} \to \Alg{S}$ are
induced from cartesian monad morphisms $S \to T$. We say a homomorphism
$f : A \to B$ of of $T$-algebras to be \emph{cartesian} if the square
\xydiagram{
  TA \ar[d]_{Th} \ar[r]^f & A \ar[d]^h \\
  TB \ar[r]_g & B
}
is a pullback square.

\begin{prop} \label{prop:cartesian-algebra-functor}
  Let $S$ and $T$ be cartesian monads.
  A monad morphism $\alpha : S \to T$ is cartesian iff the induced functor
  $\alpha^* : \Alg{T} \to \Alg{S}$ preserves cartesian homomorphisms.
\end{prop}
\begin{proof}
  Suppose $\alpha$ is cartesian and $f : A \to B$ is a cartesian homomorphism
  between algebras $\inlinearrow{TA}{\phi_A}{A}$ and
  $\inlinearrow{TB}{\phi_b}{B}$.
  Then $\alpha^*(f)$ is a cartesian homomorphism since the composite
  \xydiagram{
    SA \ar[r]^{Sf} \pullbackcorner \ar[d]_{\alpha_A} & SB \ar[d]^{\alpha_B} \\
    TA \ar[r]^{Tf} \pullbackcorner \ar[d]_{\phi_A} & TB \ar[d]^{\phi_B} \\
    A \ar[r]_f & B \\
  }
  is a pullback square.

  Conversely, suppose $\alpha^*$ preserves cartesian homomorphism. Since $T$ is
  a cartesian monad, for any
  map $g : A \to B$ the square
  \xydiagram{
    T^2 A \ar[r]^{T^2 g} \ar[d]_{\mu_A} \pullbackcorner & T^2 B \ar[d]_{\mu_A} \\
    T A \ar[r]_{T g} & T B \\
  }
  is a pullback square. Its image under $\alpha^*$ is the lower square of
  \xydiagram{
    SA \ar[r]^{Sg} \pullbackcorner \ar[d]_{S\eta_A} & SB \ar[d]^{S\eta_B} \\
    STA \ar[r]^{STg} \pullbackcorner \ar[d] & STB \ar[d] \\
    TA \ar[r]_{Tg} & TB \\
  }
  The upper square is a pullback too, so the composite square, which is the
  naturality square of $\alpha$ at $g$, is a pullback as well. Since this is
  true of any map $g$, $\alpha$ is cartesian.
\end{proof}

\subsection{Initial Algebras} \label{section:initial-algebras}

Given an endofunctor $F : \C \to \C$, an object $X \in C$ together with
an isomorphism $FX \iso X$ is called a \term{fixed point} of $F$.
The archetypal example is the natural numbers $\N$, which is a fixed point of
the endofunctor $1 + (\dash) : \Set \to \Set$. It satisfies the isomorphism
$[0, S] : 1 + \N \iso \N$ defined by $0 : 1 \to \N$ and the successor function
$S : \N \to \N$.

Such a fixed point is not unique. For example, the set
$\ExtN = \N \cup \{ \infty \}$ is also a fixed point
$[0, S] : 1 + \ExtN \iso \ExtN$, where $S(\infty) = \infty$. Although $\N$ and
$\ExtN$ are isomorphic as sets, any such isomorphism $\N \to \ExtN$ does not
commute with the isomorphisms $1 + \N \iso \N$ and $1 + \ExtN \iso \ExtN$:
any element $n \in \N$ can be written as $S(S(\ldots(0)))$ for
finitely many applications of $S$, but $\infty \in \ExtN$ cannot be written in
this way so it cannot be in the image of any map $\N \to \ExtN$ commuting
with their fixed point isomorphisms.
Therefore $\N$ and $\ExtN$ are different fixed points of $1 + (\dash)$.
However, we may instead say that $\N$ is
\term{least fixed point}, meaning there is an a unique map from $\N$ to any
other fixed point, and $\N$ is the unique fixed point having this property.

We make this precise as follows. Similar to a monad algebra, for an endofunctor
$F : \C \to \C$ and
\term{$F$-algebra} is a set $X$ and a function $FX \to X$, but it is not
required to satisfy any laws. A natural transformation $\alpha : F \to G$
similarly induces a unique functor $\alpha^* : \Alg{G} \to \Alg{F}$ lying over
$\C$. The \term{initial algebra} for an endofunctor $F$, if it exists, is the
initial object in the category of $F$-algebras.

\begin{prop}[Lambek \cite{Lambek68}]
  An initial $F$-algebra is a fixed point of $F$.
\end{prop}
\begin{proof}
  Let $\inlinearrow{FA}{a}{A}$ be the an initial $F$-algebra. Then
  $\inlinearrow{F^2A}{Fa}{FA}$ is also an $F$-algebra. By initiality, there
  exists a morphism $h : A \to FA$ such that the diagram
  \xydiagram{
    FA \ar@/_1pc/[dd]_{F1} \ar[r]^a \ar[d]^{Fh} & A \ar@/^1pc/[dd]^{1} \ar[d]_{h} \\
    FFA \ar[d]^{Fa} \ar[r]^{Fa} & FA \ar[d]_a \\
    FA \ar[r]_a & A \\
  }
  commutes. We see that $a \o h = 1$ and $h \o a = F(a \o h) = 1$, so
  $h$ is a two-sided inverse to $a$. Therefore $a$ is an isomorphism.
\end{proof}

An initial algebra does not always exist. For example, Cantor's Theorem says
that the powerset functor $\PowerSet : \Set \to \Set$ cannot have a fixed point.
The following theorem gives sufficient conditions for the existence of an
initial algebra.

\begin{thm}[Ad\'amek \cite{Adamek74}] \label{thm:initial-algebra-colimit}
  Let $\C$ be a category with an initial object, and let $F : \C \to \C$ be an
  endofunctor. If the diagram
  \labeledxydiagram{initial-algebra-colimit}{
    0 \ar[r]^{i} & F0 \ar[r]^{Fi} & F^20 \ar[r]^{F^2i} & F^30 \ar[r]^{F^3i} & \ldots
  }
  has a colimit and $F$ preserves this colimit, then this colimit carries the
  structure of an initial $F$-algebra.
\end{thm}
\begin{proof}
  Let $A$ be the colimit. The statement ``$F$ preserves this colimit'' means
  precisely that the canonical map $A \to FA$ is an isomorphism. Let
  $\alpha : FA \to A$ be its inverse; then $(A, \alpha)$ is an $F$-algebra.

  Suppose $\inlinearrow{FX}{\phi}{X}$ is another $F$-algebra. We define the legs
  $f_n : F^n0 \to B$ of a cocone over (\ref{initial-algebra-colimit}) by
  induction: the leg $f_0 : 0 \to B$ is the
  initial morphism, and the leg $f_{n + 1} : F^{n + 1}0 \to B$ is the map
  \xydiagram{
    F^{n + 1}0 \ar[r]^{F f_n} & FB \ar[r]^\beta & B.
  }
  We can verify that these legs commute with the diagram
  (\ref{initial-algebra-colimit}), so we have a unique colimit map
  $f : A \to X$. Furthermore, taking the colimit of the diagrams
  \xydiagram{
    F^{n + 1}0 \ar[r]^{Ff_n} \ar[dr]^{f_{n + 1}} & FX \ar[d]^{\phi} \\
    F^n 0 \ar[r]_{f_n} \ar[u]^{F^n i} & X
  }
  yields a diagram
  \xydiagram{
    FA \ar[r]^{Ff} & FX \ar[d]^{\phi} \\
    A \ar[r]_{f} \ar[u]^{\alpha^{-1}} & X
  }
  showing that $f$ is an $F$-algebra homomorphism.
  Since the forgetful functor $\Alg{F} \to \C$ is faithful, $\phi$ is unique
  as an $F$-algebra homomorphism as well.
\end{proof}

\subsection{Free Monads}

Given an endofunctor $F$, we may ask whether there exists a monad whose algebras
are the same as $F$-algebras. Such a monad is called a \term{free monad} for
$F$. Our main application of initial algebras is the construction of these free
monads.

Suppose $F$ is an endofunctor and $M$ is a monad. A natural
transformation $\lambda : F \to M$ induces a functor $\lambda^+$ defined as
the composite
\xydiagram{
  \MndAlg{M} \ar@{^{(}->}[r] & \EndAlg{M} \ar[r]^{\lambda^*} & \EndAlg{F}
}
lying over $\C$. Similar to Proposition \ref{prop:star-inverse}, every functor
$\MndAlg{N} \to \EndAlg{F}$ lying over $\C$ is $\lambda^+$ for some $\lambda$.

\begin{defn}[{{\cite[\S 22]{Kelly80}}}]
  Let $F : \C \to \C$ be an endofunctor. A \term{free monad} for $F$ is a monad
  $M$ and a natural transformation $\lambda : F \to M$ so that the
  functor $\lambda^+$ is an equivalence.
\end{defn}

A free monad, if it exists, satisfies the following universal property, which
justifies the ``free'' terminology.

\begin{prop}[\cite{Kelly80}] \label{prop:free-monad}
  Let $F \in \End(\C)$ be an endofunctor, and let $M$ together with
  $\lambda : F \to M$ be a free monad
  for $F$. Then for any monad $N$, composition with $\lambda$ induces a bijection
  \[
    \Mnd(\C)(M, N) \iso \End(\C)(F, N).
  \]
\end{prop}

In other words, any natural transformation $F \to N$ factors uniquely through
$\lambda$ and a monad morphism $M \to N$.

\begin{proof}
  We have a chain of isomorphisms
  \begin{align*}
    \Mnd(\C)(M, N) &\iso \Cat/\C(\MndAlg{N}, \MndAlg{M}) \\
      &\iso \Cat/\C(\MndAlg{N}, \EndAlg{F}) \\
      &\iso \End(\C)(F, N) \\
    \alpha & \mapsto \alpha^* \\
      & \mapsto \alpha^* \o \lambda^+ = (\alpha \o \lambda)^+ \\
      & \mapsto \alpha \o \lambda
  \end{align*}
  where the composite isomorphism is composition with $\lambda$.
\end{proof}

We give sufficient conditions for the existence of a free monad, using an
initial algebra construction.

\begin{thm}[{{\cite[\S 9.4]{Barr85}}}] \label{thm:free-monad-existence}
  Let $\C$ be a category with coproducts, and let $F : \C \to \C$ be an
  endofunctor. Suppose that for each object $X \in \C$, the functor
  $G_X = X + F(\dash)$ has an initial algebra. Then $F$ has a free monad.
\end{thm}
\begin{proof}
  We show that the forgetful functor $U : \Alg{F} \to \C$ is monadic. Then
  the induced monad will be the free monad for $F$, since by monadicity its
  algebras are equivalent to $F$-algebras.

  To give a left adjoint $U$, we need to find for each $X \in \C$ an
  initial object of the comma category $\comma{X}{U}$ \cite[p121]{MacLane71}.
  An object of $\comma{X}{U}$ consists of an $F$-algebra
  $\inlinearrow{FY}{f}{Y}$ and a morphism $a : X \to Y$, which is equivalently a
  $G_X$-algebra $\inlinearrow{X + FY}{[a, f]}{Y}$. Since we have an initial
  $G_X$-algebra for each $X$, we have an initial object of $\comma{X}{U}$ for
  each $X$ and thus a left adjoint $L$.

  Let $M = UL$. Explicitly, the left adjoint $L$ sends an object $X$ to the
  $F$-algebra \[ F(MX) \to X + F(MX) \iso MX, \]
  the unit $\eta_X$ of the adjunction (and the monad $M$) is the composite
  \[ X \to X + F(MX) \iso MX, \]
  and the counit $\epsilon_\phi$ at an $F$-algebra $\inlinearrow{FX}{\phi}{X}$
  is composite of the two squares
  \xydiagram{
    F(MX) \ar@/_3pc/[dd]_{LU\phi \,=\, LX} \ar[d] \ar[r]^{F\epsilon_\phi} & FX \ar@/^3pc/[dd]^{\phi} \ar[d] \\
    X + F(MX) \ar[d]^*[@]{\iso} \ar[r]^{G_X \epsilon_\phi} & X + FX \ar[d]_{[1, \phi]} \\
    MX \ar[r]_{\epsilon_\phi} & X
  }
  where $\epsilon_\phi$ is the initial map of $G_X$-algebras making the bottom
  square commute.

  Let $\lambda : FX \to MX$ be the composite
  \[ FX \xrightarrow{F\eta} F(MX) \xrightarrow{LX} MX. \]
  We show that $\lambda^+$ and the comparison functor $K^M$ are inverses,
  meaning that $M$ is monadic and that $\lambda$ has the required property.

  Let $\phi : FX \to X$ be an $F$-algebra. The equality of two composites
  $FX \to X$ in the diagram
  \xydiagram{
    FX \ar[d]_{F\eta} \ar[rd]^{1} \\
    F(MX) \ar@/_3pc/[dd]_{LU\phi} \ar[d] \ar[r]^{F\epsilon_\phi} & FX \ar@/^3pc/[dd]^{\phi} \ar[d] \\
    X + F(MX) \ar[d]^*[@]{\iso} \ar[r]^{G_X \epsilon_\phi} & X + FX \ar[d]_{[1, \phi]} \\
    MX \ar[r]_{\epsilon_\phi} & X
  }
  shows that $\lambda^+(K^M(\phi)) = \phi$, as $K^M(\phi)$ is the $T$-algebra
  $\epsilon_\phi : MX \to X$.

  Conversely, for any $M$-algebra
  $\inlinearrow{MX}{f}{X}$ the diagram
  \xydiagram{
    X + F(MX) \ar[d]^*[@]{\iso} \ar[r]^{G_X f} & X + FX \ar[d]^{[1, \lambda^+(f)]} \\
    MX \ar[r]_{f} & X
  }
  commutes, so $K^M(\lambda^+(f)) = f$ by uniqueness of the initial morphism $f$.
\end{proof}

Since colimits of functor categories are computed pointwise, if the hypotheses
of Theorem \ref{thm:free-monad-existence} are satisfied then the functor
\[
  1 + F \o (\dash) : \End(\C) \to \End(\C)
\]
has an initial algebra in $\End(\C)$, and this functor is $M$. In particular,
there is a natural isomorphism $1 + FM \iso M$, with components
$\eta : 1 \to M$ and $L : FM \to M$.

Free monads always exist for finitary endofunctors.
Let $\FinEnd(\C)$ and $\FinMnd(\C)$ be the full subcategories of
$\End(\C)$ and $\Mnd(\C)$ of finitary endofunctors and monads respectively.

\begin{prop} \label{prop:finitary-free-monad}
  Let $\C$ be a presheaf category and $F : \C \to \C$ a
  finitary endofunctor. Then the free monad on $F$ exists and is finitary.
  More generally, the forgetful functor $\FinMnd(\C) \to \FinEnd(\C)$
  has a left adjoint $\FreeMnd : \FinEnd(\C) \to \FinMnd(\C)$.
\end{prop}
\begin{proof}
  Since $F$ is finitary, each $G_X = X + F(\dash)$ is finitary also.
  The colimit (\ref{initial-algebra-colimit}) is filtered, so the colimit
  exists in $\C$ and $G_X$ preserves this colimit for each $X$. Thus the
  hypotheses of Theorem \ref{thm:free-monad-existence} are satisfied so
  $F$ has a free monad, which is the composite of the adjunction
  \[ \adjunction{\C}{\Alg{F}}{L}{U}. \]

  The left adjoint $L$ preserves colimits, so it is finitary.
  If $F$ preserves colimits of a certain shape, then the forgetful functor
  $U : \EndAlg{F} \to \C$ creates those
  colimits, so in particular if $F$ is finitary then $U$ is finitary also.
  Thus the composite $UL$ i.e.\ the free monad on $F$ is finitary.
\end{proof}

The universal property in Proposition \ref{prop:free-monad} implies that
the free monad construction
is left adjoint to the forgetful functor $\FinMnd(\C) \to \FinEnd(\C)$,
so we have an adjunction
\[ \adjunction{\FinEnd(\C)}{\FinMnd(\C)}{\FreeMnd}{U}. \]

The above propositions restrict to the case of cartesian endofunctors as well.
Let $\FinCartEnd(\C)$ and $\FinCartMnd(\C)$ be the full subcategories of
$\CartEnd(\C)$ and $\CartMnd(\C)$ of finitary objects.

\begin{thm} \label{thm:cartesian-free-monad}
  Let $\C$ be a presheaf category.
  The adjunction in Proposition \ref{prop:finitary-free-monad} restricts to
  an adjunction
  \[ \adjunction{\FinCartEnd(\C)}{\FinCartMnd(\C)}{\FreeMnd}{U}. \]
\end{thm}

\begin{proof}[Proof of Theorem \ref{thm:cartesian-free-monad}]
  We first show that if $F$ is a finitary cartesian monad, then the
  free monad $M = \FreeMnd(F)$ is cartesian as well. Suppose
  \xydiagram{
    A \pullbackcorner \ar[r] \ar[d] \pullbackcorner & B \ar[d] \\
    C \ar[r] & D
  }
  is a pullback square in $\C$. We may show by induction that for every $n$
  the square
  \xydiagram{
    G_A^n 0 \pullbackcorner \ar[r] \ar[d] \pullbackcorner & G_B^n 0 \ar[d] \\
    G_C^n 0 \ar[r] & G_D^n 0
  }
  is a pullback square,
  using the fact that $F$ preserves pullbacks and that coproducts commute
  with pullbacks.
  Taking colimits and using the commutativity of pullbacks and filtered colimits,
  we obtain a pullback square
  \xydiagram{
    MA \pullbackcorner \ar[r] \ar[d] \pullbackcorner & MB \ar[d] \\
    MC \ar[r] & MD
  }
  so the functor part of $M$ is cartesian.

  The unit $\eta : 1 \to M$ is a coproduct inclusion and so is cartesian.
  The component of the multiplication $\mu_X : M^2X \to MX$ is the initial map
  in the diagram
  \xydiagram{
    MX + F(M^2X) \ar[d]^*[@]{\iso} \ar[r]^{G_{MX}(\mu_X)} & MX + F(MX) \ar[d]^{[1, LX]} \\
    M^2X \ar[r]_{\mu_X} & X
  }
  In fact, $M^2$ is the initial algebra in $\FinCartEnd(\C)$ of
  \[
    H = M + F \o (\dash) : \FinCartEnd(\C) \to \FinCartEnd(\C).
  \]
  Since $L : FM \to M$ is a coproduct injection, it is cartesian, and thus
  so is $[1, L] : M + FM \to M$. By initiality,
  there is a unique cartesian natural transformation $\alpha$ making the diagram
  \xydiagram{
    M + FM^2 \ar[d]^*[@]{\iso} \ar[r]^{M + F\alpha} & M + FM \ar[d]^{[1, L]} \\
    M^2 \ar[r]_{\alpha} & M
  }
  commute. Since the forgetful functor $\FinCartEnd(\C) \to \FinEnd(\C)$ is
  faithful, by uniqueness $\alpha = \mu$ so $\mu$ is cartesian.

  Finally, using the characterization of cartesian monad morphisms
  in Proposition \ref{prop:cartesian-algebra-functor} in the proof of
  Proposition \ref{prop:free-monad} shows that composition
  with $\lambda$ induces a bijection
  \[
    \CartMnd(\C)(\FreeMnd(F), N) \iso \CartEnd(\C)(F, N)
  \]
  and thus $\FreeMnd : \FinCartEnd(\C) \to \FinCartMnd(\C)$ is a left adjoint.
\end{proof}

\section{Operads} \label{section:operads}

The word ``operad'' was introduced by J.P.\ May as a portmanteau of the
words ``operation'' and ``monad'' \cite{May}. Accordingly, an operad is an
algebraic theory whose elements we think of as \term{operations} of various
arity. Leinster generalizes May's operads to operads whose ``arity'' takes
values in some cartesian monad $T$, and
shows in \cite[Corollary 6.2.4]{Leinster04} that such an operad can
be characterized as a cartesian monad $P$ lying over $T$ in
the sense that we have a cartesian monad morphism $P \to T$. Recent papers
use this characterization as the definition of an operad
\cite{Garner10, Batanin11}.

The advantage of this definition, besides being extraordinarily succinct, is
that we can use the well-developed theory of monads to work with operads.
As an example, we show the free monad construction immediately yields a
construction of free operads.

To get a feel for operads, we first deconstruct the simple example of
\term{planar operads}, which are operads over the free monoid monad of Example
\ref{example:free-monoid-monad}. In Section \ref{section:globular-operads}, we
consider \term{globular operads}, which are operads over the free strict
$\omega$-category monad, and present Leinster's weak $\omega$-categories as
an algebra for such a globular operad.

\subsection{Collections and Operads}

Let $\C$ be a presheaf, and let $T$ be a finitary cartesian monad on $\C$.

\begin{defn}
  The category $\Coll{T}$ of \term{collections} over $T$ is the slice category
  $\CartEnd(\C)/T$.
\end{defn}

That is, a collection over $T$ is a cartesian endofunctor $F : \C \to \C$ with a
cartesian natural transformation $\alpha : F \natto T$. It suffices to study the
image $F1$ and the component $\alpha_1 : F1 \to T1$, as for any object $A$ we
can determine both $FA$ and $\alpha_A : FA \to TA$ via the pullback
\xydiagram{
  FA \pullbackcorner \ar[r]^{F!} \ar[d]_{\alpha_A} & F1 \ar[d]^{\alpha_1} \\
  TA \ar[r]_{T!} & T1.
}

\begin{prop}
  The functor
  \[
    \begin{array}{cccc}
      S : &\Coll{T} &\to &\C/T1 \\
        &\slicearrow{F}{\alpha}{T} &\mapsto &\slicearrow{F1}{\alpha_1}{T1}
    \end{array}
  \]
  is an equivalence.
\end{prop}
\begin{proof}
  We construct an inverse functor as sketched above. Given an object
  $\inlinearrow{P}{p}{T}$ of $\C/T1$, define the image $FA$ and the component
  $\alpha_A$ by the pullback over $T! : TA \to T1$. For $f : A \to B$, let the
  morphism $Ff$ be the universal morphism
  \xydiagram{
    FA \ar@/^1pc/[rr]^{F!} \ar@{-->}[r]_{Ff} \ar[d]_{\alpha_A} & FB \pullbackcorner
      \ar[r]_{F!} \ar[d]_{\alpha_B} & F1 \ar[d]^{\alpha_1} \\
    TA \ar@/_1pc/[rr]_{T!} \ar[r]^{Tf} & TB \ar[r]^{T!} & T1
  }
  so that the pullback lemma implies $\alpha$ is cartesian. Since
  pullbacks are unique up to isomorphism, $S$ is an equivalence.
\end{proof}

Similarly, we define an operad to be a collection whose endofunctor has a
monad structure.

\begin{defn}
  The category $\Operad{T}$ of \term{operads} over $T$ is the slice category
  $\CartMnd(\C)/T$.
  An \term{algebra} for an operad $P \natto T$ is an algebra for the
  monad $P$.
\end{defn}

\subsection{Planar Operads} \label{section:planar-operads}

An operad over the free monoid monad of Example
\ref{example:free-monoid-monad} is known as a \term{planar operad}.
Let $T$ be this free monoid monad. Since the free monoid on one generator
is the natural numbers $\N$ under addition, we have an isomorphism $T1 \iso \N$.

Let $\inlinearrow{F}{\alpha}{T}$ be a collection on $T$ and let $P = F1$. Under
the equivalence $\Coll{T} \iso \Set/T1$, the collection $\inlinearrow{F}{\alpha}{T}$ is
uniquely determined by the set $P$ and the function
$\alpha_1 : P \to \N \iso T1$. We say the set $P$ is a collection of
\term{operations}, and the function $\alpha_1 : P \to \N$ assigns
an \term{arity} to each operation. In other words, a collection over $T$ is
a finitary signature.

We write $P(n)$ for the preimage of
$\alpha_1$ at $n$, so that $P(n)$ is the set of $n$-ary operations. The
value of the functor $F$ at a set $A$ is determined by the pullback of
$\alpha_1$ along the ``length'' function $T! : TA \to \N$ as in the diagram
\xydiagram{
  FA \pullbackcorner \ar[r]^-{F!} \ar[d]_{\alpha_A} & F1 = P \ar[d]^{\alpha_1} \\
  TA \ar[r]_-{T!} & T1 \iso \N.
}
Explicitly, we have an isomorphism
\[
  FA \iso \coprod_{n \in \N} P(n) \times A^n
\]
i.e.\ an element of $FA$ is an $n$-ary operation with each input ``labeled''
by an element of $A$. We refer to the elements of $FA$ as \term{operations} as
well, and these operations have a list $TA$ as their ``arities''. Given a
list $l \in TA$ and an operation $\theta \in FA$, we say that $\theta$
\term{lies over} $l$ when $\alpha(\theta) = l$.

Suppose further that $\inlinearrow{F}{\alpha}{T}$ is an operad, i.e.\ $F$ is
a cartesian monad and $\alpha$ is a cartesian monad morphism. The unit
$\eta : 1 \to F$ is determined by the component $\eta_1 : 1 \to F1 = P$,
i.e.\ by an element of $P$. The commutativity of the diagram
\xydiagram{
  & F1 = P \ar[dd]^{\alpha_1} \\
  1 \ar[ur]^{\eta_1} \ar[dr]_{\eta_1} \\
  & T1 \iso \N
}
implies that this element is a unary operation. We call this
operation the \term{identity} operation.
The multiplication $\mu : F^2 \to F$
is determined by the component $\mu_1 : F^2 1 \to 1$ i.e.\ $\mu_1 : FP \to P$.
Using the commutativity of
\xydiagram{
  FP \ar[r]^{\mu_1} \ar[d]_{\alpha\alpha_1} & P \ar[d]^{\alpha_1} \\
  T^2 1 \iso T\N \ar[r]_{\mu_1} & T1 \iso \N
}
we can show that $\mu_1$ takes an $n$-ary operation labeled with ``inner'' operations
of arity $(k_1, \ldots, k_n)$ to an operation of arity $k_1 + \ldots + k_n$.
Thus $\mu$ specifies a way to compose operations. The monad laws impose
associativity and identity laws (e.g.\ composing an operation $\theta \in P(n)$
with the identity operation results in $\theta$), although it is tedious to
write these laws explicitly.

An algebra for $\inlinearrow{F}{\alpha}{T}$ consists of a set $A$ and a
function $f : FA \to A$ that \term{applies} an operation to its labels to
yield an element of $A$. Equivalently, this means for every unlabeled operation
$\theta \in P$ of arity $n$, we have a function
\[
  \theta_A : A^n \to A
\]
and these functions should respect the composition structure of the operad.

\subsection{Free Operads}

There is a faithful functor $U : \Operad{T} \to \Coll{T}$ that forgets the monad
structure (equivalently, takes the component at $1$).

\begin{prop}
  The functor $U : \Operad{T} \to \Coll{T}$ has a left adjoint
  $\FreeOperad : \Coll{T} \to \Operad{T}$.
\end{prop}

This adjoint is simply the free monad construction.

\begin{proof}
  Suppose we have a collection $\inlinearrow{F}{\alpha}{T}$.
  Since $T$ is finitary, by Proposition \ref{prop:cartesian-list} the functor
  $F$ is finitary also. By Theorem \ref{thm:cartesian-free-monad}, there
  exists a finitary and cartesian free monad $\FreeMnd(F)$ of $F$.
  The natural transformation $\alpha$ induces a cartesian monad morphism
  $\inlinearrow{\FreeMnd(F)}{\adjunct{\alpha}}{T}$ and thus an operad.
  Any map of collections $F \to P$ to an operad $P$ must factor uniquely through
  the map
  \xydiagram{
    F \ar[dr]_{\alpha} \ar[rr]^{\lambda_F} & & \FreeMnd(F) \ar[dl]^{\adjunct{\alpha}} \\
    & T
  }
  so $\FreeOperad : \Coll{T} \to \Operad{T}$ is a left adjoint with $\lambda_F$
  as its unit.
\end{proof}

\begin{example}
  Suppose $T$ is the free monoid monad on $\Set$. We can define a collection
  \[ \slicearrow{1}{2}{\N \iso T1} \]
  which picks out $2 \in \N$. The induced collection is the cartesian functor
  $P = (\dash)^2 : \Set \to \Set$, where the operad map sends an ordered pair $(a, b)$ to the
  two-element list $(a, b)$. The
  free operad on this collection is the set
  of (finite) binary trees with leaves labeled by elements of $A$,
  where the map $\FreeMnd(P) \to TA$ collects the labels into a list in order.
  An algebra for this free operad is equivalently a
  $(\dash)^2$-algebra, i.e.\ merely a set $X$ and a function
  $X \times X \to X$ that is not required to satisfy any laws (this structure
  is known as a \term{magma}).

  More generally, suppose $\Sigma \to \N$ is a finitary signature. Then an algebra
  for the free operad on this collection is a $\Sigma$-algebra in the sense
  of Example \ref{example:equational-presentation}.
\end{example}

\section{Globular Operads} \label{section:globular-operads}

The goal of this section is to present Leinster's definition of a weak
$\omega$-category. Its components are \term{globular sets}, which
form the cellular structure of weak $\omega$-categories;
\term{globular operads}, which encode the operations we can
perform in weak $\omega$-categories; and \term{contractions}, which specifies
that a globular operad is suitably weak.

\subsection{Globular Sets} \label{section:globular-sets}

In Example \ref{example:free-category-monad}, we defined a category as an
algebra for a monad on the category of graphs. An $\omega$-category should be
an algebra for a monad on the category of ``$\omega$-graphs'', or
\term{globular sets}. A globular set consists of a set of 0-cells (objects),
1-cells (morphisms) between 0-cells, 2-cells (transformations) between 1-cells,
and so on, so that each $n$-cell has one source $(n - 1)$-cell and one
target $(n - 1)$-cell. We say $n$-cells are \term{parallel} if they have
the same source and target. An $n$-cell's source and target must be parallel;
this is analogous to the condition that functors $F$ and $G$ must have the same
domain and codomain categories in order for a natural transformation $F \to G$ to be
well-defined.

As with the category of graphs, we may define the category of globular sets as a
presheaf category.

\begin{defn}
  The \term{globe category} $\G$ is the category whose objects are the
  natural numbers $\N = \{0, 1, \ldots\}$ and whose morphisms generated by
  \xydiagram{
    0 \ar@<0.5ex>[r]^{\sigma_0} \ar@<-0.5ex>[r]_{\tau_0} &
    1 \ar@<0.5ex>[r]^{\sigma_1} \ar@<-0.5ex>[r]_{\tau_1} &
    2 \ar@<0.5ex>[r]^{\sigma_2} \ar@<-0.5ex>[r]_{\tau_2} &
    \ldots
  }
  such that for all $n \in \N$
  \begin{align*}
    \sigma_{n + 1} \o \sigma_n &= \tau_{n + 1} \o \sigma_n \\
    \sigma_{n + 1} \o \tau_n &= \tau_{n + 1} \o \tau_n.
  \end{align*}
  The category $\GSet$ of \term{globular sets} is the presheaf category
  $[\op{\G}, \Set]$.
\end{defn}

This definition agrees with the informal description of globular sets above.
Let $X : \op{\G} \to \Set$ be a globular set. For each $n \in \N$, the image
$X(n)$ is the set of $n$-cells. We call the images of $\sigma_n$ and $\tau_n$
\[ s_n, t_n : X(n + 1) \to X(n) \]
respectively, and they map an $(n + 1)$-cell to its source and target $n$-cells.
For any $n$-cell $c$, the equations
\begin{align*}
  s_{n - 2}(s_{n - 1}(c)) &= s_{n - 2}(t_{n - 1}(c)) \\
  t_{n - 2}(s_{n - 1}(c)) &= t_{n - 2}(t_{n - 1}(c))
\end{align*}
hold, so the source and target of $c$ are parallel.

Each topological space $S$ has an associated globular set of homotopies, its
\term{globular nerve}.
Let $\disk{(\dash)} : \G \to \Top$ be the functor that sends the object $n$ to
the $n$-dimensional disk $\disk{n}$, and $\sigma_n$ and $\tau_n$ to
the embedding of $\disk{n}$ as the upper or lower hemisphere of the boundary of
$\disk{n + 1}$. The globular nerve functor $N : \Top \to \GSet$ sends a
space $S$ to the globular set
\[
  NS = \Top(\disk{(\dash)}, S)
\]
so that an $n$-cell is a map $\disk{n} \to S$.
We can describe the low dimensional cells of $NS$ as follows:
\begin{itemize}
\item a 0-cell is a point,
\item a 1-cell is a path between a source point and a target point,
\item a 2-cell is a disk i.e.\ a homotopy between two paths relative to their endpoints,
\item a 3-cell is a ball i.e.\ a homotopy between two hemisphere disks relative to
  their intersection at the equator,
\end{itemize}
and so on.

The globular nerve $N$ has a right adjoint $\realize{\dash} : \GSet \to \Top$,
the \term{geometric realization} \cite[p240]{Leinster04}. The geometric
realization $|X|$ of a globular set $X$ is a CW complex constructed by making an
$n$-disk for each $n$-cell, and gluing the boundary of each $n$-disk to its
source and target $(n - 1)$-disks.

The image of Yoneda embedding
$\yoneda : \G \to \GSet$ at $n$ is the
\term{standard $n$-globe}\footnote{%
The notation $G_n$ follows Brown \cite{Brown08}; alternatively, Street \cite{Street00} and
Garner \cite{Garner10} use $\mathbf{2}_n$.}
\[ G_n = \yoneda(n), \]
whose geometric realization is the disk $\disk{n}$. By abuse of notation, we
give the injections $\yoneda(\sigma_n), \yoneda(\tau_n)$ the names
$\sigma_n, \tau_n : G_n \toto G_{n + 1}$. As a globular set, $G_n$ consists
of a single $n$-cells, two $(n - 1)$-cells
(the source and target of the $n$-cell),
two $(n - 2)$-cells (the source and target of the $(n - 1)$-cells), and so on.
By the Yoneda Lemma, for any globular set $X$ we have a bijection
$\GSet(G_n, X) \iso X(n)$, natural in both $n$ and $X$.

We construct a sequence of globular sets $\del_n$ whose geometric
realization is the $(n - 1)$-sphere\footnote{%
The notation $\del_n$ is due to Garner \cite{Garner10}, who constructs these
globular sets in an equivalent but slightly more complex way.}, and
maps $\boundary_n : \del_n \to G_n$ whose realization embeds the
boundary $\sphere{(n - 1)}$ into $\disk{n}$. For example, the 2-sphere is the
union of the upper and lower hemispheres of the 3-ball, glued along their
common boundary 1-sphere. This motivates the following definition.

\begin{defn}
We construct the sequence of globular sets $\del_n$ and maps
$\boundary_n : \del_n \to G_n$ by induction. Let $\del_0$
be the initial object and $i_0 : \del_0 \to G_0$ be the initial morphism. For $n
\in \N$, let $\del_{n + 1}$ and $\boundary_{n + 1}$ be the pushout
and universal morphism in
\labeledxydiagram{eq:boundary-pushout}{
  \del_n \ar[d]_{\boundary_n} \ar[r]^{\boundary_n} & G_n \ar[d] \ar@/^/[rdd]^{\tau_n} \\
  G_n \ar[r] \ar@/_/[rrd]_{\sigma_n}
    & \del_{n + 1} \pushoutcorner \ar@{.>}[rd]|{\boundary_{n + 1}} \\
  & & G_{n + 1}.
}
\end{defn}

As a globular set, $\del_n$ consists of two $(n - 1)$-cells, two $(n - 2)$-cells,
and so on.
Since the realization of $\del_n$ is topologically a pair of ``parallel''
$(n - 1)$-disks glued along their boundary, we expect that a map of globular
sets $\del_n \to X$ consists of a pair of parallel $(n - 1)$-cells of $X$. Let
\[
  \mathrm{Par}(X, n) = \begin{cases}
    \singleton & n = 0 \\
    \{ (a, b) \in X(n - 1) \times X(n - 1) \mid
      \text{$a$ and $b$ are parallel} \} & n > 0
  \end{cases}
\]
where $(n + 1)$ cells $c$ and $d$ are parallel if $(s(c), t(c)) = (s(d), t(d))$,
and all 0-cells are parallel. Since the source and target of any cell are
parallel, the function
$d : X(n) \to \mathrm{Par}(X, n)$ that sends an $(n + 1)$-cell to its
source and target $(s(c), t(c))$ is well-defined.

\begin{prop} \label{prop:boundary-representable}
  Let $X$ be a globular set. We have a natural bijection
  \[
    \GSet(\del_n, X) \iso \mathrm{Par}(X, n).
  \]
  such that the diagram
  \xydiagram{
    \GSet(G_n, X) \ar[d]_{\GSet(i_n, X)} \ar[r]^{\iso} & X(n) \ar[d]^d \\
    \GSet(\del_n, X) \ar[r]^{\iso} & \mathrm{Par}(X, n)
  }
  commutes in $\Set$, where the top arrow is the natural bijection from the
  Yoneda Lemma.
\end{prop}
\begin{proof}
  The sets
  $\GSet(\del_0, X)$ and $\mathrm{Par}(X, 0)$ both have one element, so the
  desired bijection is trivial.

  The hom-functor
  $\GSet(\dash, X)$ turns colimits in $\GSet$ into limits in $\Set$, so applying
  this functor on diagram \eqref{eq:boundary-pushout} gives us a diagram
  \labeledxydiagram{eq:boundary-pushout-image}{
    \GSet(G_{n + 1}, X) \ar@{.>}[rd] \ar@/_/[rdd] \ar@/^/[rrd] \\
    & \GSet(\del_{n + 1}, X) \ar[r] \ar[d] \pullbackcorner &
      \GSet(G_n, X) \ar[d] \\
    & \GSet(G_n, X) \ar[r] & \GSet(\del_n, X)
  }

  Note that set-theoretic definition of $\mathrm{Par}(X, n)$ above is a pullback in
  $\Set$ of $X(n - 1)$ with itself, so the following
  diagram commutes in $\Set$.
  \labeledxydiagram{eq:set-pullback}{
    X(n + 1) \ar@{.>}[rd]|{d} \ar@/_/[rdd]_s \ar@/^/[rrd]^t \\
    & \mathrm{Par}(X, n + 1) \ar[r] \ar[d] \pullbackcorner &
      X(n) \ar[d]^{d} \\
    & X(n) \ar[r]_{d} & \mathrm{Par}(X, n)
  }

  Using the Yoneda Lemma, the inductive hypothesis, and the universal property
  of pullbacks, we can construct bijections between the corresponding objects
  of \eqref{eq:boundary-pushout-image} and \eqref{eq:set-pullback} so that
  everything in sight commutes. This yields the desired bijection
  \[
    \GSet(\del_{n + 1}, X) \iso \mathrm{Par}(X, n + 1).
  \]
  and the bijections at the dashed arrow yields
  \xydiagram{
    \GSet(G_{n + 1}, X) \ar[d]_{\GSet(i_{n + 1}, X)} \ar[r]^{\iso} & X(n + 1) \ar[d]^d \\
    \GSet(\del_{n + 1}, X) \ar[r]^{\iso} & \mathrm{Par}(X, n + 1)
  }
  as required.
\end{proof}

\subsection{Globular Operads}

Although they are perhaps less interesting,
strict $\omega$-categories can be given the following explicit (if verbose)
definition, which we have slightly modified from \cite[p. 22]{Leinster04}.

\begin{defn} \label{defn:strict-category}
A \term{strict $\omega$-category} is a globular set $A$ equipped with the
following operations.
\begin{itemize}
  \item For all $0 \le k < n$, a function
    $\fo_k : A(n) \times_{A(p)} A(n) \to A(n)$ on the pullback
    \xydiagram{
      \cdot \pullbackcorner \ar[r] \ar[d] & A(n) \ar[d]^{s_k} \\
      A(n) \ar[r]_{t_k} & A(k)
    }
    that maps \term{$k$-composable} $n$-cells $a, b \in A(n)$ such that
    $t_k(a) = s_k(b)$ to their \term{$k$-composite} $n$-cell $a \fo_k b$.

  \item For all $n$, a function $i : A(n) \to A(n + 1)$ that maps
    an $n$-cell $a$ to its \term{identity} $(n + 1)$-cell $i(a)$.
    We write $i_k : A(k) \to A(n)$ for the repeated application $i^{n - k}$.
\end{itemize}

These operations must satisfy axioms the following axioms.
\begin{enumerate}[label=(\alph*)]
  \item (boundaries of composites) for compatible $n$-cells $a, b \in A(n)$
  \begin{align*}
    s_l(a \fo_k b) &= s_l(a) &\text{ and } t_l(a \fo_k b) &= t_l(b) & \text{ if } k = l \\
    s_l(a \fo_k b) &= s_l(a) \fo_k s_l(b) &\text{ and } t_l(a \fo_k b) &= t_l(a) \fo_k t_l(b) & \text{ if } k < l.
  \end{align*}

  \item (boundaries of identities) For $a \in A(k)$ then
  \[ s_k(i_k(a)) = a = t_k(i_k(a)). \]

  \item (associativity) For compatible $n$-cells $a, b, c \in A(n)$ then
  \[ a \fo_k (b \fo_k c) = (a \fo_k b) \fo_k c. \]

  \item (identities) For $0 \le k < n$ and $x \in A(n)$ then
  \[ i_k(t_k(x)) \fo_k x = x = x \fo_k i_k(s_k(x)). \]

  \item (binary exchange) For $0 \le k < l < n$ and compatible $n$-cells
  $a, b, c, d \in A(n)$ then
  \[ (a \fo_k b) \fo_l (c \fo_k d) = (a \fo_l c) \fo_k (b \fo_l d). \]

  \item (nullary exchange) For $0 \le k < n$ and compatible $n$-cells
  $a, b \in A(n)$ then
  \[ i(a) \fo_k i(b) = i(a \fo_k b). \]
\end{enumerate}

\end{defn}

If $A$ and $B$ are strict $\omega$-categories then a
\term{strict $\omega$-functor} is a map $f : A \to B$ of globular sets that
commutes with the composition and identity operations.
We may form an (ordinary) category $\StrCat$ consisting of strict
$\omega$-categories and strict $\omega$-functors.

Let $U : \StrCat \to \GSet$ be the forgetful functor that extracts the
underlying globular set of a strict $\omega$-category. Leinster shows in
\cite[Appendix F]{Leinster04} that this functor $U$ is monadic, meaning in particular
that it has a left adjoint. This left adjoint sends a globular set $A$ to
the \term{free strict $\omega$-category} on $A$. This adjunction induces a monad
$T$, such that $TA$ is the underlying globular set of the free strict
$\omega$-category on $A$.

For $n \in \N$, the elements of $TA(n)$ are called
\term{$n$-dimensional pasting diagrams}, or simply $n$-pasting diagrams. These
represent formal composites of elements of $A$, and they will play the
role of ``arities'' for globular operations similar to the role the free
monoid monad plays for planar operads in Section \ref{section:planar-operads}.

Leinster gives a formal description of pasting diagrams in
\cite[\S 8.1]{Leinster04} and \cite[Appendix F]{Leinster04}. We briefly give
an intuitive description of the monad $T$ here. An element of $TA(2)$ might
look like
\[
  \begin{tikzcd}
    A \arrow[r, bend left=75, looseness=2, "f", ""{below, name=f}]
      \arrow[r, bend left=35, "g"{above, near start}, ""{above, name=g1}, ""{below, name=g2}]
      \arrow[r, bend right=35, "h"{below, near start}, ""{above, name=h1}, ""{below, name=h2}]
      \arrow[r, bend right=75, looseness=2, "i"{below}, ""{above, name=i}]
    & B
      \arrow[r, "j"]
    & C
      \arrow[r, bend left, "k", ""{below, name=k}]
      \arrow[r, bend right, "l"{below}, ""{above, name=ell}]
    & D
    \arrow[Rightarrow, from=f, to=g1, "\alpha"]
    \arrow[Rightarrow, from=g2, to=h1, "\beta"]
    \arrow[Rightarrow, from=h2, to=i, "\gamma"]
    \arrow[Rightarrow, from=k, to=ell, "\delta"]
    & \in TA(2)
  \end{tikzcd}
\]
given 2-cells $\alpha, \beta, \gamma, \delta \in A(2)$, 1-cells
$f, g, h, \ldots \in A(1)$, and 0-cells $A, B, C, D \in A(0)$ such that
that the sources and targets of the cells in the pasting diagram
match, i.e.
\begin{align*}
  s(f) &= A & s(\alpha) &= f \\
  t(f) &= B & t(\alpha) &= g \\
  s(g) &= A & s(\beta) &= g \\
  t(g) &= B & t(\beta) &= h \\
  s(h) &= A & s(\gamma) &= h \\
  t(h) &= B & t(\gamma) &= i \\
  \vdots & & \vdots
\end{align*}
and so on. The pasting diagram represents the formal composite 2-cell
\[
  (\alpha \fo_1 \beta \fo_1 \gamma) \fo_0 \id{j} \fo_0 \gamma.
\]
The source and target of the above pasting diagram are the
1-pasting diagrams
\[
  \begin{tikzcd}
    A \arrow[r, "f"]
    & B
      \arrow[r, "j"]
    & C
      \arrow[r, "k"]
    & D
    & \in TA(1)
  \end{tikzcd}
\]
and
\[
  \begin{tikzcd}
    A \arrow[r, "i"]
    & B
      \arrow[r, "j"]
    & C
      \arrow[r, "l"]
    & D
    & \in TA(1)
  \end{tikzcd}
\]
respectively.
Note that cells is the diagram are determined by other the source or target of
other cells
(e.g. the top left 1-cell must be labeled with the source of $\alpha$, namely $f$).
For brevity, we will often omit unnecessary labels when drawing pasting diagrams.

The $\omega$-category structure on $TA$ is given by ``pasting'' two diagrams
together along a common boundary. For example, we have composites
\begin{align*}
  \left(
  \begin{tikzcd}[ampersand replacement=\&]
    A
      \arrow[r, bend left=50, "f"{above}, ""{below, name=f1}]
      \arrow[r, bend right=50, "g"{below}, ""{above, name=f2}]
    \& B
      \arrow[r, bend left=50, "i"{above}, ""{below, name=g1}]
      \arrow[r, bend right=50, "j"{below}, ""{above, name=g2}]
    \& C
    \arrow[Rightarrow, from=f1, to=f2, "\alpha"]
    \arrow[Rightarrow, from=g1, to=g2, "\gamma"]
  \end{tikzcd}
  \right) \fo_1 \left(
  \begin{tikzcd}[ampersand replacement=\&]
    A
      \arrow[r, bend left=50, "g"{above}, ""{below, name=f1}]
      \arrow[r, bend right=50, "h"{below}, ""{above, name=f2}]
    \& B
      \arrow[r, "j"{above}, ""{below, name=g1}]
    \& C
    \arrow[Rightarrow, from=f1, to=f2, "\beta"]
  \end{tikzcd}
  \right) = \left(
  \begin{tikzcd}[ampersand replacement=\&]
    A
      \arrow[r, bend left=75, "f"{above}, ""{below, name=f1}]
      \arrow[r, "g"{near start}, ""{above, name=f2}, ""{below, name=g1}]
      \arrow[r, bend right=75, "h"{below}, ""{above, name=g2}]
    \& B
      \arrow[r, bend left=50, "i"{above}, ""{below, name=h1}]
      \arrow[r, bend right=50, "j"{below}, ""{above, name=h2}]
    \& C
    \arrow[Rightarrow, from=f1, to=f2, "\alpha"]
    \arrow[Rightarrow, from=g1, to=g2, "\beta"]
    \arrow[Rightarrow, from=h1, to=h2, "\gamma"]
  \end{tikzcd}
  \right) \\
  \left(
  \begin{tikzcd}[ampersand replacement=\&]
    A \arrow[r, bend left=75, looseness=2, "f", ""{below, name=f}]
      \arrow[r, bend left=35, "g"{above, near start}, ""{above, name=g1}, ""{below, name=g2}]
      \arrow[r, bend right=35, "h"{below, near start}, ""{above, name=h1}, ""{below, name=h2}]
      \arrow[r, bend right=75, looseness=2, "i"{below}, ""{above, name=i}]
    \& B
      \arrow[r, "j"]
    \& C
    \arrow[Rightarrow, from=f, to=g1, "\alpha"]
    \arrow[Rightarrow, from=g2, to=h1, "\beta"]
    \arrow[Rightarrow, from=h2, to=i, "\gamma"]
  \end{tikzcd}
  \right) \fo_0 \left(
  \begin{tikzcd}[ampersand replacement=\&]
    C
      \arrow[r, bend left, "k", ""{below, name=k}]
      \arrow[r, bend right, "l"{below}, ""{above, name=ell}]
    \& D
    \arrow[Rightarrow, from=k, to=ell, "\delta"]
  \end{tikzcd}
  \right) = \left(
  \begin{tikzcd}[ampersand replacement=\&]
    A \arrow[r, bend left=75, looseness=2, "f", ""{below, name=f}]
      \arrow[r, bend left=35, "g"{above, near start}, ""{above, name=g1}, ""{below, name=g2}]
      \arrow[r, bend right=35, "h"{below, near start}, ""{above, name=h1}, ""{below, name=h2}]
      \arrow[r, bend right=75, looseness=2, "i"{below}, ""{above, name=i}]
    \& B
      \arrow[r, "j"]
    \& C
      \arrow[r, bend left, "k", ""{below, name=k}]
      \arrow[r, bend right, "l"{below}, ""{above, name=ell}]
    \& D
    \arrow[Rightarrow, from=f, to=g1, "\alpha"]
    \arrow[Rightarrow, from=g2, to=h1, "\beta"]
    \arrow[Rightarrow, from=h2, to=i, "\gamma"]
    \arrow[Rightarrow, from=k, to=ell, "\delta"]
  \end{tikzcd}
  \right).
\end{align*}
It may help to stack the pasting diagrams vertically when composing with
$\fo_1$; we have written them horizontally for clarity.

Suppose $\pi$ is the 1-pasting diagram
\[
  \begin{tikzcd}
    A \arrow[r, "f"]
    & B
      \arrow[r, "g"]
    & C
      \arrow[r, "h"]
    & D
    & \in TA(1)
  \end{tikzcd}
\]
i.e.\ the formal composite $\pi = f \fo_0 g \fo_0 h$. The identity pasting
diagram $\id{\pi}$ is the ``degenerate'' 2-pasting diagram
\[
  \begin{tikzcd}
    A \arrow[r, "f"]
    & B
      \arrow[r, "g"]
    & C
      \arrow[r, "h"]
    & D
    & \in TA(2)
  \end{tikzcd}
\]
which represents the formal composite
\[ \id{\pi} = \id{f \fo_0 g \fo_0 h} = \id{f} \fo_0 \id{g} \fo_0 \id{h}. \]
Note that $\id{\pi}$ has the same pictorial representation as $\pi$. When drawing
the diagrams, we will disambiguate the dimension by specifying which set the
pasting diagram is an element of, e.g.\ with ``$\in TA(1)$'' or
``$\in TA(2)$'' as above.

Finally, we describe the monad structure of $T$. The unit $\eta : 1 \to T$
takes a cell to a pasting diagram consisting of only that cell, so
e.g.\ if $\alpha \in A(2)$ is a 2-cell, then $\eta(\alpha)$ is the 2-pasting
diagram
\[
  \begin{tikzcd}[ampersand replacement=\&]
    \pdpt \arrow[r, bend left=60, "s(\alpha)", ""{below, name=f}]
      \arrow[r, bend right=60, "t(\alpha)"{below}, ""{above, name=g}] \& \pdpt.
    \arrow[Rightarrow, from=f, to=g, "\alpha"]
  \end{tikzcd}
\]
The multiplication $\mu : TT \to T$ takes a pasting diagram of pasting diagrams,
such as
\[
  \begin{tikzcd}
    \pdpt
      \arrow[r, bend left=75, ""{below, name=f1}]
      \arrow[r, ""{above, name=f2}, ""{below, name=g1}]
      \arrow[r, bend right=75, ""{above, name=g2}]
    & \pdpt
      \arrow[r, bend left, ""{below, name=k}]
      \arrow[r, bend right, ""{above, name=ell}]
    & \pdpt
    \arrow[Rightarrow, from=f1, to=f2, "\pi_1"]
    \arrow[Rightarrow, from=g1, to=g2, "\pi_2"]
    \arrow[Rightarrow, from=k, to=ell, "\pi_3"]
    & \in TTA(2)
  \end{tikzcd}
\]
where
\begin{align*}
  \pi_1 &=
  \begin{tikzcd}[ampersand replacement=\&]
    A \arrow[r, bend left=75, "f", ""{below, name=f}]
      \arrow[r, "g"{above, near start}, ""{above, name=g1}, ""{below, name=g2}]
      \arrow[r, bend right=75, "h"{below, near start}, ""{above, name=h1}, ""{below, name=h2}]
    \& B
    \arrow[Rightarrow, from=f, to=g1, "\alpha"]
    \arrow[Rightarrow, from=g2, to=h1, "\beta"]
  \end{tikzcd}
  & \in TA(2) \\
  \pi_2 &=
  \begin{tikzcd}[ampersand replacement=\&]
    A \arrow[r, bend left, "h"{above}, ""{below, name=h2}]
      \arrow[r, bend right, "i"{below}, ""{above, name=i}]
    \& B
    \arrow[Rightarrow, from=h2, to=i, "\gamma"]
  \end{tikzcd}
  & \in TA(2) \\
  \pi_3 &=
  \begin{tikzcd}[ampersand replacement=\&]
    B
      \arrow[r, "j"]
    \& C
      \arrow[r, bend left, "k", ""{below, name=k}]
      \arrow[r, bend right, "l"{below}, ""{above, name=ell}]
    \& D
    \arrow[Rightarrow, from=k, to=ell, "\delta"]
  \end{tikzcd}
  & \in TA(2)
\end{align*}
and composes them as the outer pasting diagram indicates. Thus the pasting diagram
above is mapped to the composite $(\pi_1 \fo_1 \pi_2) \fo_0 \pi_3$ which is
the pasting diagram
\[
  \begin{tikzcd}
    A \arrow[r, bend left=75, looseness=2, "f", ""{below, name=f}]
      \arrow[r, bend left=35, "g"{above, near start}, ""{above, name=g1}, ""{below, name=g2}]
      \arrow[r, bend right=35, "h"{below, near start}, ""{above, name=h1}, ""{below, name=h2}]
      \arrow[r, bend right=75, looseness=2, "i"{below}, ""{above, name=i}]
    & B
      \arrow[r, "j"]
    & C
      \arrow[r, bend left, "k", ""{below, name=k}]
      \arrow[r, bend right, "l"{below}, ""{above, name=ell}]
    & D
    \arrow[Rightarrow, from=f, to=g1, "\alpha"]
    \arrow[Rightarrow, from=g2, to=h1, "\beta"]
    \arrow[Rightarrow, from=h2, to=i, "\gamma"]
    \arrow[Rightarrow, from=k, to=ell, "\delta"]
    & \in TA(2).
  \end{tikzcd}
\]

Leinster additional proves in \cite[Appendix F]{Leinster04} that the
free strict $\omega$-category monad $T$ is finitary and cartesian. Thus we may
speak of operads over $T$, which are known as \term{globular operads}.

We use similar terminology as Section \ref{section:planar-operads} to describe
the structure of globular operads. Let $\inlinearrow{P}{\alpha}{T}$ be a
globular operad. For a globular set $A$ and $n \in N$, we call the elements of
$PA(n)$ \term{$n$-dimensional operations} or $n$-operations. Each $n$-operation
$\theta \in PA(n)$ has an ``arity'' $n$-pasting diagram $\alpha(\theta) \in
TA(n)$; in this case we say that $\theta$ \term{lies over} the pasting diagram
$\alpha(\theta)$. An algebra for a globular operad consists of a way to
\emph{apply} each $n$-dimensional operation to form a result $n$-cell in $A$.

\subsection{Contractions}

Note that $T$ itself is naturally a globular operad (indeed,
the terminal globular operad) whose algebras are strict $\omega$-categories.
A weak $\omega$-category should be an algebra for an operad that is ``weakly
equivalent'' to $T$ in some sense. We take some inspiration from topology in
order to define what this means.

A map $q : X \to Y$ of topological spaces is \term{weak homotopy equivalence}
\cite{Jardine87} iff for each commuting square
\xydiagram{
  S^{n - 1} \ar[d]_i \ar[r] & X \ar[d]^q \\
  \disk{n} \ar[r] & Y
}
there exists a map $\sigma : \disk{n} \to X$ that fills the diagonal of the square
and makes the two resulting triangles commute up to homotopy:
\xydiagram{
  S^{n - 1} \ar[d]_i \ar[r] & X \ar[d]^q \\
  \disk{n} \ar@{.>}[ur]^\sigma \ar[r] & Y.
}
A weak homotopy equivalence very roughly means that the fibers of the map $q$
are contractible \cite{Smale57}.
This characterization of weak homotopy equivalences is reminiscent of the
Kan lifting property of simplicial homotopy theory \cite{SHT99}. Thinking of a map
$\disk{n} \to X$ as an $n$-dimensional homotopy, given a homotopy $\disk{n} \to Y$
with a corresponding ``shape'' $S^{n - 1} \to X$, we can ``lift'' to a homotopy
$\disk{n} \to X$ that fills the shape.
We can extend this notion of weak equivalence to globular operads, using
the globular set $\del_n$ instead of the sphere $S^{n - 1}$, and the
standard $n$-globe $G_n$ instead of the disk $\disk{n}$.

Let $\inlinearrow{P}{\alpha}{T}$ be a globular collection. Suppose we have a
pasting diagram $\pi \in TA(n)$ and parallel operations
$\theta^-, \theta^+ \in PA(n - 1)$ such that $\alpha(\theta^-) = s(\pi)$
and $\alpha(\theta^+) = t(\pi)$. This defines a diagram
\xydiagram{
  \del_n \ar[r]^{(\theta^-, \theta^+)} \ar[d]_{i_n} & PA \ar[d]^{\alpha_A} \\
  \G_n \ar[r]_{\pi} & TA \\
}
in $\GSet$. A \term{contraction} on $\inlinearrow{P}{\alpha}{T}$ allows us to
lift $\pi$ into an operation $\lambda \in PA(n)$ such that
$s(\lambda) = \theta^-$, $t(\lambda) = \theta^+$, and $\alpha(\lambda) = \pi$.
\[
  \begin{tikzcd}
    \pdpt \arrow[r, bend left=60, "\theta^-", ""{below, name=f}]
      \arrow[r, bend right=60, "\theta^+"{below}, ""{above, name=g}] & \pdpt
    \arrow[Rightarrow, from=f, to=g, "\lambda"]
  \end{tikzcd}
  \overset{\alpha}{\mapsto}
  \begin{tikzcd}
    \pdpt \arrow[r, bend left=60, "s(\pi)", ""{below, name=f}]
      \arrow[r, bend right=60, "t(\pi)"{below}, ""{above, name=g}] & \pdpt
    \arrow[Rightarrow, from=f, to=g, "\pi"]
  \end{tikzcd}
\]
This is the same as giving a diagonal filler
$\lambda : G_n \to PA$ that makes the two triangles commute in the diagram
\xydiagram{
  \del_n \ar[r]^{(\theta^-, \theta^+)} \ar[d]_{i_n} & PA \ar[d]^{\alpha_A} \\
  \G_n \ar@{.>}[ur]^{\lambda} \ar[r]_{\pi} & TA.
}

Thus we have the following definition, which is equivalent to the definition
of contraction that Leinster gives in \cite[\S 9.1]{Leinster04}.

\begin{defn}
  Let $\inlinearrow{P}{\alpha}{T}$ be a globular collection. A
  \emph{contraction} on $P$ is a function that for each $n \ge 1$ assigns each
  commuting square
  \xydiagram{
    \del_n \ar[r]^{(\theta^-, \theta^+)} \ar[d]_{i_n} & PA \ar[d]^{\alpha_A} \\
    \G_n \ar[r]_{\pi} & TA \\
  }
  to an $n$-operation $\lambda : G_n \to PA$ making the diagram
  \xydiagram{
    \del_n \ar[r]^{(\theta^-, \theta^+)} \ar[d]_{i_n} & PA \ar[d]^{\alpha_A} \\
    \G_n \ar@{.>}[ur]^{\lambda} \ar[r]_{\pi} & TA
  }
  commute.
\end{defn}

An \term{operad-with-contraction} is an operad together with a contraction on its
underlying collection. A map of operads-with-contraction is an operad map
\xydiagram{
  P \ar[dr]_{\alpha} \ar[rr]^{\phi} & & Q \ar[dl]^{\beta} \\
  & T
}
that
commutes with the contractions, in the sense that the lifts of the left and
outer squares commute in any diagram of the form
\xydiagram{
  \del_n \ar[r]^{(\theta^-, \theta^+)} \ar[d]_{i_n} & PA \ar[d]^(0.7){\alpha_A} \ar[r]^{\phi_A} & QA \ar[d]^{\beta_A} \\
  \G_n \ar@{.>}[ur]^{\lambda_P} \ar@{.>}[urr]_(0.7){\lambda_Q} \ar[r]_{\pi} & TA \ar@{=}[r] & TA.
}

Leinster shows abstractly in \cite[Appendix G]{Leinster04} that there exists an
initial operad-with-contraction $L$; Cheng \cite{Cheng10} gives a
``dimension-by-dimension'' construction of this initial operad-with-contraction,
which proceed by alternately adding in contraction lifts and operadic
compositions at each dimension. This operad is universal among
operads-with-contraction in the sense that an algebra for any
operad-with-contraction is also an algebra for $L$. Thus, we present Leinster's
definition of a weak $\omega$-category.

\begin{defn}
  A weak $\omega$-category is an algebra for the initial operad-with-contraction
  $L$.
\end{defn}

\section{Weak $\omega$-Category Operations} \label{section:weak-operations}

It is not immediately obvious why the definition of weak $\omega$-category in
the previous section is indeed a suitable definition. Thus we demonstrate how
weak $\omega$-category structure arises from this definition. In particular, we
construct composites of morphisms, associativity and coherence isomorphisms, and
an Eckmann-Hilton braiding.

First we introduce some notation. For a globular set $X$ and $a, b \in X(n)$
are $n$-cells, we write $f : a \to b$ to indicate that $f \in X(n + 1)$ is
an $(n + 1)$-cell with source $s(f) = a$ and target $t(f) = b$. We will use
this notation both for cells in a weak $\omega$-category and for operations
in a globular operad.

Let $\apply{\dash} : LA \to A$ be an $L$-algebra i.e.\ a weak $\omega$-category.
That is, if $\theta : \alpha \to \beta \in LA(n)$ is a globular operation,
we write $\apply{\theta} : \apply{\alpha} \to \apply{\beta}$ for its composite
in $A(n)$. All of our work in this section will be in the context of this
weak $\omega$-category. Generally, we will use Latin letters (e.g.\ $a$, $B$, $f$)
to denote cells in the in $\omega$-category $A$, and Greek letters
(e.g.\ $\alpha$, $\theta$, $\Phi$) for operations in $LA$.
As a special case, we let $\eta : 1 \to L$ and $\mu : L^2 \to L$ stand for
the monad operations.

\subsection{Composites}

We show how the contraction structure allows us to ``lift'' strict
$\omega$-category operations (namely, composition $\dash \fo_k \dash$ and
identity $\id{\dash}$) into weak $\omega$-category operations.
We start with composition of 1-cells $\fo_0$ along their source and target
0-cells.

This is perhaps an involved way to define composition, given that most
definitions of higher categories (including Definition
\ref{defn:strict-category}) begin by stipulating composition directly. However,
defining composition in this way gives us easy proofs of coherence laws, as
we see in Section \ref{section:coherence}.

First, we need $\omega$-category composites of operations (not to be confused
with operadic composition i.e.\ multiplication with $\mu$). Let
$\alpha, \beta, \gamma \in LA(0)$ be 0-operations, and suppose we have
1-operations $\theta : \alpha \to \beta$ and $\phi : \beta \to \gamma$. We can
form the pasting diagram of operations
\[
  \begin{tikzcd}
    \alpha \ar{r}{\theta} & \beta \ar{r}{\phi} & \gamma & \in TLA(1)
  \end{tikzcd}
\]
as $\pi = \eta_T \theta \fo_0 \eta_T \phi : \eta_T \alpha \to \eta_T \gamma$,
where $\eta_T : 1 \to T$ is the unit of the monad $T$. Now $\eta \alpha$ and
$\eta \gamma$ are parallel 0-cells in $LLA$ that lie over
$s(\pi)$ and $t(\pi)$ respectively, so we may form the contraction lift
\xydiagram{
  \del_1 \ar[d] \ar[r]^{(\eta \alpha, \eta \gamma)} & LLA \ar[d] \\
  G_1 \ar[r]_{\pi} \ar@{.>}[ru]^{\psi} & TLA
}
and thus we have a 1-cell $\psi : \eta \alpha \to \eta \gamma \in LLA(1)$.
Using the multiplication, we obtain a 1-cell
$\mu \psi : \alpha \to \gamma \in LA(1)$. We write this $\mu \psi$
as $\theta \fo_0 \phi : \alpha \to \gamma$.

Now, suppose $a, b, c \in A(0)$ and $f : a \to b$ and $g : b \to c$. Then we
have an operation
\[
  \eta f \fo_0 \eta g : \eta a \to \eta c
\]
and applying this operation yields a 1-cell
\[
  \apply{\eta f \fo_0 \eta g} : a \to c
\]
in $A(1)$, since $\apply{\eta a} = a$ and $\apply{\eta c} = c$ by the
monad algebra laws. This is the composite 1-cell of $f$ and $g$.
From this point forward, we will suppress the use of $\eta$ so e.g.\ we write
the composite 1-cell of $f$ and $g$ as $\apply{f \fo_0 g}$.

Now suppose we additional have a 0-cell $d \in A(0)$ and $h : b \to c$. We
can form the operation
\[
  (f \fo_0 g) \fo_0 h : a \to d \in LA(1).
\]
We stress that this is a \emph{different} operation than the operation
\[
  \apply{f \fo_0 g} \fo_0 h : a \to d \in LA(1).
\]
The former operation $(f \fo_0 g) \fo_0$ lies over the pasting diagram
\[
  \begin{tikzcd}
    a \ar{r}{f} & b \ar{r}{g} & c\ar{r}{h} & d & \in TA(1)
  \end{tikzcd}
\]
while the latter operation $\apply{f \fo_0 g} \fo_0 h$ lies over the pasting
diagram
\[
  \begin{tikzcd}
    a \ar{r}{\apply{f \fo_0 g}} & c\ar{r}{h} & d & \in TA(1).
  \end{tikzcd}
\]
However, the applications of these 1-cells are equal in $A(1)$, meaning
\begin{equation} \label{eq:L-algebra-law}
  \apply{(f \fo_0 g) \fo_0 h} = \apply{\apply{f \fo_0 g} \fo_0 h},
\end{equation}
which follows from the monad algebra law
\xydiagram{
  LLA \ar[d]_{\mu_A} \ar[r]^{L\apply{\dash}} & LA \ar[d]^{\apply{\dash}} \\
  LA \ar[r]_{\apply{\dash}} & A.
}
More generally, we can ``remove'' nested applications. This fact is
crucial to relate operations lying over different pasting diagrams, and
plays an important role in Section \ref{section:eckmann-hilton}.

In much the same way, we may lift any strict $\omega$-category operation to a
weak $\omega$-category operation. As one more example, if we have operations
\begin{align*}
  \alpha, \beta, \gamma & \in LA(0) \\
  \theta_1, \theta_2 &: \alpha \to \beta \\
  \phi_1, \phi_2 &: \beta \to \gamma \\
  \Theta &: \theta_1 \to \theta_2 \\
  \Phi &: \phi_1 \to \phi_2
\end{align*}
we may form the operation
\[ \Theta \fo_0 \Phi : \theta_1 \fo_0 \phi_1 \to \theta_2 \fo_0 \phi_2 \]
by lifting the pasting diagram
\[
  \begin{tikzcd}
    \mathclap{\alpha}
      \arrow[r, bend left=50, "\theta_1"{above}, ""{below, name=f1}]
      \arrow[r, bend right=50, "\theta_2"{below}, ""{above, name=f2}]
    & \mathclap{\beta}
      \arrow[r, bend left=50, "\phi_1"{above}, ""{below, name=g1}]
      \arrow[r, bend right=50, "\phi_2"{below}, ""{above, name=g2}]
    & \mathclap{\gamma}
    \arrow[Rightarrow, from=f1, to=f2, "\Theta"]
    \arrow[Rightarrow, from=g1, to=g2, "\Phi"]
    & {\in TLA(2),}
  \end{tikzcd}
\]
and we may form $\id{\alpha} : \alpha \to \alpha$ by lifting the identity
1-cell $\eta_T \alpha \in TLA(1)$. This implies that if we have
$\omega$-category cells
\begin{align*}
  a, b, c & \in LA(0) \\
  f_1, f_2 &: a \to b \\
  g_1, g_2 &: b \to c \\
  F &: f_1 \to f_2 \\
  G &: g_1 \to g_2
\end{align*}
we can form $\omega$-category cells
\[
  \apply{F \fo_0 G} : \apply{f_1 \fo_0 g_1} \to \apply{f_2 \fo_0 g_2}
\]
and
\[
  \apply{\id{a}} : a \to a.
\]

\subsection{Associativity, Units, and Exchange}

Our next task is to show that the strict $\omega$-category laws, namely
associativity, units, and exchange, hold up to \term{equivalence}.
Intuitively, an equivalence should mean a ``maximally weak isomorphism'', but
a precise definition is tricky to pin down. The idea is a cell $f : a \to b$
is an equivalence if there exists a cell $f^{-1} : b \to a$ and equivalences
\begin{align*}
  f \fo f^{-1} &\to \apply{id(a)} \\
  f^{-1} \fo f &\to \apply{id(b)}.
\end{align*}
However, this ``definition'' is circular, as we need ``equivalences all the way
up''. Nevertheless, we will still be able to offer evidence that the cells we
construct are indeed equivalences.

In this section, we construct the following cells:
\begin{itemize}
  \item An \term{associator} $A_{f,g,h} : \apply{(f \fo_0 g) \fo_0 h} \to \apply{f \fo_0 (g \fo_0 h)}$,
  \item A \term{left unitor} $L_f : \apply{\id{a} \fo_0 f} \to f$,
  \item A \term{right unitor} $R_f : \apply{f \fo_0 \id{b}} \to f$, and
  \item An \term{exchanger} \[
    E_{F_1, F_2, G_1, G_2}
    : \apply{(F_1 \fo_0 G_1) \fo_1 (F_2 \fo_0 G_2)}
    \to \apply{(F_1 \fo_1 F_2) \fo_0 (G_1 \fo_1 G_2)}.
    \]
\end{itemize}
Of course, there should be
equivalences like these for all dimensions, but the
construction of these cells will immediately generalize to other dimensions.

These laws all come from lifting identity pasting diagrams. For example,
note that $(f \fo_0 g) \fo_0 h$ and $f \fo_0 (g \fo_0 h)$ both lie over
the same pasting diagram
\[
  \begin{tikzcd}
    a \ar{r}{f} & b \ar{r}{g} & c \ar{r}{h} & d & \in TA(1).
  \end{tikzcd}
\]
Thus, let $\pi$ the identity 2-pasting diagram on the 1-pasting diagram above,
i.e.\ the degenerate 2-pasting diagram
\[
  \begin{tikzcd}
    a \ar{r}{f} & b \ar{r}{g} & c \ar{r}{h} & d & \in TA(2).
  \end{tikzcd}
\]
Then $(f \fo_0 g) \fo_0 h$ lies over $s(\pi)$ and $f \fo_0 (g \fo_0 h)$ lies
over $t(\pi)$, so we may lift $\pi$ to an operation
\[
  \alpha_{f,g,h} : (f \fo_0 g) \fo_0 h \to f \fo_0 (g \fo_0 h)
\]
and applying this operation yields a 2-cell
\[
  A_{f,g,h} = \apply{\alpha_{f,g,h}} : \apply{(f \fo_0 g) \fo_0 h} \to \apply{f \fo_0 (g \fo_0 h)}.
\]

Since $s(\pi) = t(\pi)$, we can switch the source and target so that $\pi$ to an operation
\[
  \alpha_{f,g,h}^{-1} : f \fo_0 (g \fo_0 h) \to (f \fo_0 g) \fo_0 h
\]
and applying this operation yields a 2-cell
\[
  A_{f,g,h}^{-1} = \apply{\alpha_{f,g,h}^{-1}} : \apply{f \fo_0 (g \fo_0 h)} \to \apply{(f \fo_0 g) \fo_0 h}.
\]
We use the name $A_{f,g,h}^{-1}$ to suggest that it is an inverse to $A_{f,g,h}$.
Indeed, note
\[
  \apply{A_{f,g,h} \fo_1 A_{f,g,h}^{-1}}
    = \apply{\apply{\alpha_{f,g,h}} \fo_1 \apply{\alpha_{f,g,h}^{-1}}}
    = \apply{\alpha_{f,g,h} \fo_1 \alpha_{f,g,h}^{-1}}.
\]
Since $\alpha_{f,g,h} \fo_1 \alpha_{f,g,h}^{-1}$ and $\id{(f \fo_0 g) \fo_0 h}$
both lie over the identity 2-pasting diagram
\[
  \begin{tikzcd}
    a \ar{r}{f} & b \ar{r}{g} & c \ar{r}{h} & d & \in TA(2),
  \end{tikzcd}
\]
lifting the identity 3-pasting diagram into a 3-operation
\[
  \alpha_{f,g,h} \fo_1 \alpha_{f,g,h}^{-1} \to \id{(f \fo_0 g) \fo_0 h}
\]
whose application is a 3-cell
\[
  \apply{A_{f,g,h} \fo_1 A_{f,g,h}^{-1}} \to \apply{\id{(f \fo_0 g) \fo_0 h}}.
\]
Similarly, we can construct a 3-cell
\[
  \apply{A_{f,g,h}^{-1} \fo_1 A_{f,g,h}} \to \apply{\id{f \fo_0 (g \fo_0 h)}},
\]
and since these 3-cells are both lifts of identity pasting diagrams, we may
further show that these have inverses, and so on. Thus $A_{f,g,h}$ is an
equivalence.

In general, whenever two parallel operations $\theta$ and $\phi$ lie over
the same pasting diagram $\pi$, we may lift $\id{\pi}$ to an equivalence
$\apply{\theta} \to \apply{\phi}$ (or $\apply{\phi} \to \apply{\theta}$). Thus, since the operations
$\id{a} \fo_0 f$, $f$, and $f \fo_0 \id{a}$ all lie over the pasting diagram
\[
  \begin{tikzcd}
    a \ar{r}{f} & b & \in TA(1),
  \end{tikzcd}
\]
we have operations
\begin{align*}
  \lambda_f &: \id{a} \fo_0 f \to f \\
  \rho_f &: f \fo_0 \id{b} \to f.
\end{align*}
and thus equivalence 2-cells
\begin{align*}
  L_f &= \apply{\lambda_f} : \apply{\id{a} \fo_0 f} \to f \\
  R_f &= \apply{\rho_f} : \apply{f \fo_0 \id{b}} \to f.
\end{align*}
Similarly, since $(F_1 \fo_0 G_1) \fo_1 (F_2 \fo_0 G_2)$ and
$(F_1 \fo_1 F_2) \fo_0 (G_1 \fo_1 G_2)$ lie over the same pasting diagram
\[
  \begin{tikzcd}
    \pdpt
      \arrow[r, bend left=75, ""{below, name=f1}]
      \arrow[r, ""{above, name=f2}, ""{below, name=g1}]
      \arrow[r, bend right=75, ""{above, name=g2}]
    & \pdpt
      \arrow[r, bend left=75, ""{below, name=h1}]
      \arrow[r, ""{above, name=h2}, ""{below, name=j1}]
      \arrow[r, bend right=75, ""{above, name=j2}]
    & \pdpt
    \arrow[Rightarrow, from=f1, to=f2, "F_1"]
    \arrow[Rightarrow, from=g1, to=g2, "F_2"]
    \arrow[Rightarrow, from=h1, to=h2, "G_1"]
    \arrow[Rightarrow, from=j1, to=j2, "G_2"]
  \end{tikzcd}
\]
we have an equivalence 3-cell
\[
  E_{F_1, F_2, G_1, G_2}
  : \apply{(F_1 \fo_0 G_1) \fo_1 (F_2 \fo_0 G_2)}
  \to \apply{(F_1 \fo_1 F_2) \fo_0 (G_1 \fo_1 G_2)}.
\]

\subsection{Coherence} \label{section:coherence}

We demonstrate that the equivalences we constructed in the previous section
are \emph{coherent}, which in some sense means that all diagrams built out of
these equivalences commute. The most famous of these diagrams are the
\term{pentagon diagram}
\[
  \begin{tikzpicture}
    \node (P0) at (90:3.0cm) {$\apply{(f \fo_0 g) \fo_0 (h \fo_0 i)}$};
    \node (P1) at (90+72:3.0cm) {$\apply{((f \fo_0 g) \fo_0 h) \fo_0 i}$};
    \node (P2) at (90+2*72:3.0cm) {$\apply{(f \fo_0 (g \fo_0 h)) \fo_0 i}$};
    \node (P3) at (90+3*72:3.0cm) {$\apply{f \fo_0 ((g \fo_0 h) \fo_0 i)}$};
    \node (P4) at (90+4*72:3.0cm) {$\apply{f \fo_0 (g \fo_0 (h \fo_0 i))}$};
    \draw
    (P1) edge[->] node[left, near end, outer sep=8pt] {$A_{f,g,\apply{h \fo_0 i}}$} (P0)
    (P0) edge[->] node[right, near start, outer sep=8pt] {$A_{\apply{f \fo_0 g},h,i}$} (P4)
    (P1) edge[->] node[left] {$\apply{\id{f} \fo_0 A_{g,h,i}}$} (P2)
    (P2) edge[->] node[below, outer sep=6pt] {$A_{f, \apply{g \fo_0 h}, i}$} (P3)
    (P3) edge[->] node[right] {$\apply{A_{f,g,h} \fo_0 \id{i}}$} (P4);
  \end{tikzpicture}
\]
and the \term{triangle diagram}
\[
\begin{tikzcd}
  \apply{(f \fo_0 \id{b}) \fo_0 g}
    \ar[swap]{rd}{\apply{R_f \fo_0 \id{g}}}
    \ar{rr}{A_{f,\apply{\id{b}},g}}
  & & \apply{f \fo_0 (\id{b} \fo_0 g)} \ar{ld}{\apply{\id{f} \fo_0 L_f}}  \\
  & \apply{f \fo_0 g},
\end{tikzcd}
\]
which are traditionally part of the definition of a monoidal category
\cite{MacLane71}.
We show that these diagrams commute up to equivalence, meaning we exhibit
3-cells from the composite of one path in the diagrams to the other path.
Specifically, we need 3-cells
\begin{align*}
  \apply{A_{f,g,\apply{h \fo_0 i}} \fo_1 A_{\apply{f \fo_0 g},h,i}}
    \to
    \apply{\apply{\id{f} \fo_0 A_{g,h,i}} \fo_1 (A_{f, \apply{g \fo_0 h}, i} \fo_1 \apply{A_{f,g,h} \fo_0 \id{i}})}
\end{align*}
and
\begin{align*}
  \apply{A_{f,\apply{\id{b}},g} \fo_1 \apply{\id{f} \fo_0 L_f}} \to \apply{R_f \fo_0 \id{g}}.
\end{align*}

With $\alpha$, $\lambda$, and $\rho$ as in the previous section, removing nested applications results in 3-cells
\begin{align*}
  \apply{\alpha_{f,g,h \fo_0 i} \fo_1 \alpha_{f \fo_0 g,h,i}}
    \to
    \apply{(\id{f} \fo_0 \alpha_{g,h,i}) \fo_1 (\alpha_{f, g \fo_0 h, i} \fo_1 (\alpha_{f,g,h} \fo_0 \id{i}))}
\end{align*}
and
\begin{align*}
  \apply{\alpha_{f,\id{b},g} \fo_1 (\id{f} \fo_0 \lambda_f}) \to \apply{\rho_f \fo_0 \id{g}}.
\end{align*}
Note that both operations
\[
  \alpha_{f,g,h \fo_0 i} \fo_1 \alpha_{f \fo_0 g,h,i}
  \; \text{ and } \;
  (\id{f} \fo_0 \alpha_{g,h,i}) \fo_1 (\alpha_{f, g \fo_0 h, i} \fo_1 (\alpha_{f,g,h} \fo_0 \id{i}))
\]
lie over the identity 2-pasting diagram
\[
  \begin{tikzcd}
    a \ar{r}{f} & b \ar{r}{g} & c \ar{r}{h} & d \ar{r}{i} & e & \in TA(2)
  \end{tikzcd}
\]
so by the remarks in the previous section we have an equivalence between
their applications. Similarly, both
\[
  \alpha_{f,\id{b},g} \fo_1 (\id{f} \fo_0 \lambda_f)
  \; \text{ and } \;
  \apply{\rho_f \fo_0 \id{g}}
\]
lie over the identity 2-pasting diagram
\[
  \begin{tikzcd}
    a \ar{r}{f} & b \ar{r}{g} & c & \in TA(2)
  \end{tikzcd}
\]
so we have an equivalence between their applications as well.
Thus the pentagon and triangle diagram commute up to equivalence.

In general, proving coherence laws simply require lifting an appropriate
identity pasting diagram. This approach also us to show \emph{naturality}
diagrams commute as well. For example, suppose $F : f \to f'$ is a 2-cell
in $A$. We expect that the diagram
\[
  \begin{tikzcd}[column sep=1.5in]
    \apply{(f \fo_0 g) \fo_0 h}
        \ar{r}{\apply{(F \fo_0 \id{g}) \fo_0 \id{h}}}
        \ar{d}{A_{f,g,h}}
      & \apply{(f' \fo_0 g) \fo_0 h} \ar{d}{A_{f',g,h}} \\
    \apply{f \fo_0 (g \fo_0 h)} \ar{r}[swap]{\apply{F \fo_0 (\id{g} \fo_0 \id{h})}}
      & \apply{f' \fo_0 (g \fo_0 h)}
  \end{tikzcd}
\]
should commute, as the associator should be a natural equivalence. Indeed, both
the top right and bottom left paths are applications of operations that
lie over the pasting diagram
\[
  \begin{tikzcd}
    a
      \arrow[r, bend left=50, "f"{above}, ""{below, name=f1}]
      \arrow[r, bend right=50, "f'"{below}, ""{above, name=f2}]
    & \mathclap{b} \arrow[r, "g"]
    & c \arrow[r, "h"]
    & d
    \arrow[Rightarrow, from=f1, to=f2, "F"]
    & \in TA(2).
  \end{tikzcd}
\]
so we can lift this into an equivalence making the diagram commute.

Thus, by using the same contraction structure to define both composition and
associators, unitors, etc., the coherence laws come essentially for free.

\subsection{The Eckmann-Hilton Argument} \label{section:eckmann-hilton}

In its original form, the Eckmann and Hilton showed that if a set is
equipped with two monoid structures such that one is a homomorphism for the
other, then in fact the two monoid structures coincide that the resulting monoid
is commutative \cite{Eckmann62}. In the context of higher category theory,
the monoid structures are composition of endomorphisms, and they are related
by the exchange law \cite{BaezDolan95}.

Let $\point \in A(0)$ be a 0-cell. We consider the endomorphisms of the
1-cell $\apply{\id{\point}} : \point \to \point$, i.e.\ 2-cells
$\apply{\id{\point}} \to \apply{\id{\point}}$. Given any two such 2-cells
$a, b : \apply{\id{\point}} \to \apply{\id{\point}}$, we may form
the 1-composite $\apply{a \fo_1 b}$ by composing along $\apply{\id{\point}}$,
but we may also form the 0-composite
$\apply{a \fo_0 b}$ by composing along $\point$.

The composition is weakly commutative in the sense that there exists a 2-cell
\[
  \braid{a,b} : \apply{a \fo_1 b} \to \apply{b \fo_1 a},
\]
which we will now construct. The idea is the use the extra dimension to
``rotate'' the two 2-cells around each other using the identity and
exchange laws.

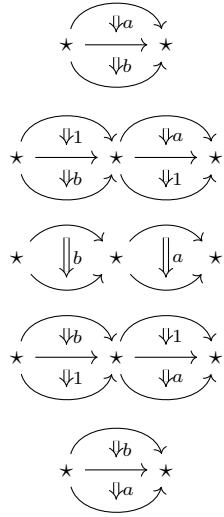
\[
\begin{gathered}
  \begin{tikzcd}[ampersand replacement=\&]
    \point
      \arrow[r, bend left=75, ""{below, name=f1}]
      \arrow[r, ""{above, name=f2}, ""{below, name=g1}]
      \arrow[r, bend right=75, ""{above, name=g2}] \& \point
    \arrow[Rightarrow, from=f1, to=f2, "a"]
    \arrow[Rightarrow, from=g1, to=g2, "b"]
  \end{tikzcd} \\
  \begin{tikzcd}[ampersand replacement=\&]
    \point
      \arrow[r, bend left=75, ""{below, name=f1}]
      \arrow[r, ""{above, name=f2}, ""{below, name=g1}]
      \arrow[r, bend right=75, ""{above, name=g2}]
    \& \point
      \arrow[r, bend left=75, ""{below, name=h1}]
      \arrow[r, ""{above, name=h2}, ""{below, name=j1}]
      \arrow[r, bend right=75, ""{above, name=j2}]
    \& \point
    \arrow[Rightarrow, from=f1, to=f2, "1"]
    \arrow[Rightarrow, from=g1, to=g2, "b"]
    \arrow[Rightarrow, from=h1, to=h2, "a"]
    \arrow[Rightarrow, from=j1, to=j2, "1"]
  \end{tikzcd} \\
  \begin{tikzcd}[ampersand replacement=\&]
    \point
      \arrow[r, bend left=50, ""{below, name=f1}]
      \arrow[r, bend right=50, ""{above, name=f2}]
    \& \point
      \arrow[r, bend left=50, ""{below, name=g1}]
      \arrow[r, bend right=50, ""{above, name=g2}]
    \& \point
    \arrow[Rightarrow, from=f1, to=f2, "b"]
    \arrow[Rightarrow, from=g1, to=g2, "a"]
  \end{tikzcd} \\
  \begin{tikzcd}[ampersand replacement=\&]
    \point
      \arrow[r, bend left=75, ""{below, name=f1}]
      \arrow[r, ""{above, name=f2}, ""{below, name=g1}]
      \arrow[r, bend right=75, ""{above, name=g2}]
    \& \point
      \arrow[r, bend left=75, ""{below, name=h1}]
      \arrow[r, ""{above, name=h2}, ""{below, name=j1}]
      \arrow[r, bend right=75, ""{above, name=j2}]
    \& \point
    \arrow[Rightarrow, from=f1, to=f2, "b"]
    \arrow[Rightarrow, from=g1, to=g2, "1"]
    \arrow[Rightarrow, from=h1, to=h2, "1"]
    \arrow[Rightarrow, from=j1, to=j2, "a"]
  \end{tikzcd} \\
  \begin{tikzcd}[ampersand replacement=\&]
    \point
      \arrow[r, bend left=75, ""{below, name=f1}]
      \arrow[r, ""{above, name=f2}, ""{below, name=g1}]
      \arrow[r, bend right=75, ""{above, name=g2}] \& \point
    \arrow[Rightarrow, from=f1, to=f2, "b"]
    \arrow[Rightarrow, from=g1, to=g2, "a"]
  \end{tikzcd}
\end{gathered}
\]

We decompose the construction of $\apply{a \fo_1 b} \to \apply{b \fo_0 a}$
into three steps. First, we form a weak 3-operation
\[
  \theta_1 : a \fo_1 b \to
    \Big(\id{\apply{\id{\point}}} \fo_0 a \Big) \fo_1
    \Big(b \fo_0 \id{\apply{\id{\point}}} \Big)
\]
from lifting the identity 3-pasting diagram on the 2-pasting diagram
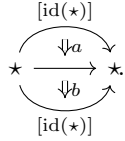
\[
  \begin{tikzcd}[ampersand replacement=\&]
    \point
      \arrow[r, bend left=75, "\apply{\id{\point}}"{above}, ""{below, name=f1}]
      \arrow[r, ""{above, name=f2}, ""{below, name=g1}]
      \arrow[r, bend right=75, "\apply{\id{\point}}"{below}, ""{above, name=g2}] \& \point\mathclap{.}
    \arrow[Rightarrow, from=f1, to=f2, "a"]
    \arrow[Rightarrow, from=g1, to=g2, "b"]
  \end{tikzcd}
\]
Second, we form an exchange operation
\[
  \theta_2 :
    \Big(\apply{\id{\id{\point}}} \fo_0 a \Big) \fo_1
    \Big(b \fo_0 \apply{\id{\id{\point}}} \Big)
    \to
    \Big(\apply{\id{\id{\point}}} \fo_1 b \Big) \fo_0
    \Big(a \fo_1 \apply{\id{\id{\point}}} \Big)
\]
by lifting the identity 3-pasting diagram on the following 2-pasting diagram,
where we say $1 = \apply{\id{\id{\point}}}$ for brevity.
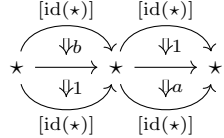
\[
  \begin{tikzcd}[ampersand replacement=\&]
    \point
      \arrow[r, bend left=75, "\apply{\id{\point}}"{above}, ""{below, name=f1}]
      \arrow[r, ""{above, name=f2}, ""{below, name=g1}]
      \arrow[r, bend right=75, "\apply{\id{\point}}"{below}, ""{above, name=g2}]
    \& \point
      \arrow[r, bend left=75, "\apply{\id{\point}}"{above}, ""{below, name=h1}]
      \arrow[r, ""{above, name=h2}, ""{below, name=j1}]
      \arrow[r, bend right=75, "\apply{\id{\point}}"{below}, ""{above, name=j2}]
    \& \point
    \arrow[Rightarrow, from=f1, to=f2, "b"]
    \arrow[Rightarrow, from=g1, to=g2, "1"]
    \arrow[Rightarrow, from=h1, to=h2, "1"]
    \arrow[Rightarrow, from=j1, to=j2, "a"]
  \end{tikzcd}
\]
Lastly, we form the weak 3-composite
\[
  \theta_3 :
    \Big(\id{\id{\point}} \fo_1 b \Big) \fo_0
    \Big(a \fo_1 \id{\id{\point}} \Big)
    \to
    b \fo_0 a
\]
by lifting the identity 3-pasting diagram on the 2-pasting diagram
\[
  \begin{tikzcd}[ampersand replacement=\&]
    \point
      \arrow[r, bend left=50, ""{below, name=f1}]
      \arrow[r, bend right=50, ""{above, name=f2}]
    \& \point
      \arrow[r, bend left=50, ""{below, name=g1}]
      \arrow[r, bend right=50, ""{above, name=g2}]
    \& \point.
    \arrow[Rightarrow, from=f1, to=f2, "b"]
    \arrow[Rightarrow, from=g1, to=g2, "a"]
  \end{tikzcd}
\]
Note that the operations
\begin{align*}
  \theta_1 &: a \fo_1 b \to
    \Big(\id{\apply{\id{\point}}} \fo_0 a \Big) \fo_1
    \Big(b \fo_0 \id{\apply{\id{\point}}} \Big) \\
  \theta_2 &:
    \Big(\apply{\id{\id{\point}}} \fo_0 a \Big) \fo_1
    \Big(b \fo_0 \apply{\id{\id{\point}}} \Big)
    \to
    \Big(\apply{\id{\id{\point}}} \fo_1 b \Big) \fo_0
    \Big(a \fo_1 \apply{\id{\id{\point}}} \Big) \\
  \theta_3 &:
    \Big(\id{\id{\point}} \fo_1 b \Big) \fo_0
    \Big(a \fo_1 \id{\id{\point}} \Big)
    \to
    b \fo_0 a
\end{align*}
are \emph{not}
composable. Their sources and targets have different bracketings of
$\id{\id{\point}}$, depending on if we want to think of the identity as
a cell in $A$ or as a weak operation in $LA$. However, similar to the remarks
for equation (\ref{eq:L-algebra-law}), the $L$-algebra laws imply that
\begin{align*}
  \apply{\theta_1} &: \apply{a \fo_1 b} \to
    \Big[ \Big(\id{\id{\point}} \fo_0 a \Big) \fo_1
    \Big(b \fo_0 \id{\id{\point}} \Big) \Big ] \\
  \apply{\theta_2} &:
    \Big[ \Big(\id{\id{\point}} \fo_0 a \Big) \fo_1
    \Big(b \fo_0 \id{\id{\point}} \Big) \Big ]
    \to
    \Big[ \Big(\id{\id{\point}} \fo_1 b \Big) \fo_0
    \Big(a \fo_1 \id{\id{\point}} \Big) \Big ] \\
  \apply{\theta_3} &:
    \Big[ \Big(\id{\id{\point}} \fo_1 b \Big) \fo_0
    \Big(a \fo_1 \id{\id{\point}} \Big) \Big ]
    \to
    \apply{b \fo_0 a}
\end{align*}
are in fact composable. Thus we may form the composite 3-cell
\[
  \apply{\apply{\theta_1} \fo_2 \apply{\theta_2} \fo_2 \apply{\theta_3} }
    : \apply{a \fo_1 b} \to \apply{b \fo_0 a}.
\]
In an entirely similar manner, we can construct 3-cells
\begin{align*}
  \theta_4 &:
    b \fo_0 a
    \to
    \Big(b \fo_1 \id{\id{\point}}\Big) \fo_0
    \Big(\id{\id{\point}} \fo_1 a \Big) \\
  \theta_5 &:
    \Big(b \fo_1 \apply{\id{\id{\point}}} \Big) \fo_0
    \Big(\apply{\id{\id{\point}}} \fo_1 a \Big)
    \to
    \Big(b \fo_0 \apply{\id{\id{\point}}} \Big) \fo_1
    \Big(\apply{\id{\id{\point}}} \fo_0 a \Big) \\
  \theta_6 &:
    \Big(b \fo_0 \id{\apply{\id{\point}}} \Big) \fo_1
    \Big(\id{\apply{\id{\point}} \fo_0 a} \Big)
    \to
    b \fo_1 a
\end{align*}
so that
\[
  \braid{a, b} =
    \apply{\apply{\theta_1} \fo_2 \apply{\theta_2} \fo_2 \apply{\theta_3} \fo_2
      \apply{\theta_4} \fo_2 \apply{\theta_5} \fo_2 \apply{\theta_6}}
      : \apply{a \fo_1 b} \to \apply{b \fo_1 a}.
\]
as required.

All the operations $\theta_i$ are all lifts of identity pasting
diagrams, so we may form $\theta_1^{-1}$,
$\theta_2^{-2}$, etc.\ by switching the source and target of the lifts. This
yields a 3-cell
\[
  \braid{a, b}^{-1} =
    \apply{\apply{\theta_6^{-1}} \fo_2 \apply{\theta_5^{-1}} \fo_2 \apply{\theta_4^{-1}} \fo_2
      \apply{\theta_3^{-1}} \fo_2 \apply{\theta_2^{-1}} \fo_2 \apply{\theta_1^{-1}}}
      : \apply{b \fo_1 a} \to \apply{a \fo_1 b}
\]
that is inverse to $\braid{a,b}$. Note that switching $a$ and $b$ gives us
a 3-cell
\[
  \braid{b, a}^{-1} : \apply{a \fo_1 b} \to \apply{b \fo_1 a}
\]
with the same source and target as $\braid{a, b}$. However, $\braid{a, b}$
and $\braid{b, a}^{-1}$ are not equivalent in general. They correspond to
different ``paths'' around the circle
\[
  \begin{tikzpicture}
    \draw (0,0) circle (2cm);
    \node [fill=white, inner sep=5pt, align=left] at (-2,0) {$a \fo_0 b$};
    \node [fill=white, inner sep=5pt, align=left] at (2,0) {$b \fo_0 a$};
    \node [fill=white, inner sep=5pt, align=left] at (0,2) {$a$ \\ $\fo_1$ \\ $b$};
    \node [fill=white, inner sep=5pt, align=left] at (0,-2) {$b$ \\ $\fo_1$ \\ $a$};
    \draw[thick, ->] ([shift={(0,0)}]75:2.8) arc (75:-75:2.8) node[midway, right]{$\braid{a,b}$};
    \draw[thick, ->] ([shift={(0,0)}]105:2.8) arc (105:255:2.8) node[midway, left]{$\braid{b,a}^{-1}$};
  \end{tikzpicture}
\]
as $\braid{a, b}$ rotates $a$ clockwise around $b$, while
$\braid{b, a}^{-1}$ rotates $a$ counterclockwise around $b$.

Pictorially, we can represent $\braid{a, b}$ as rotating two ``strings''
around in 3-dimensional space, forming a 3-dimensional ``string diagram''
(cf.\ \cite{Joyal91})
\begin{equation} \label{eq:braid}
  \begin{tikzpicture}
    \draw (0,0) .. controls (1,0) and (1,1) .. (2,1);
    \draw[-,line width=10pt,draw=white] (0,1) .. controls (1,1) and (1,0) .. (2,0);
    \draw (0,1) .. controls (1,1) and (1,0) .. (2,0);
    \draw (0,1) node[left] {$a$};
    \draw (0,0) node[left] {$b$};
    \draw (2,1) node[right] {$b$};
    \draw (2,0) node[right] {$a$};
  \end{tikzpicture}
\end{equation}
The picture for $\braid{b, a}^{-1}$ is formed by reversing $\braid{a, b}$,
and swapping the two strings:
\begin{equation} \label{eq:braid-inverse}
  \begin{tikzpicture}
    \draw (0,1) .. controls (1,1) and (1,0) .. (2,0);
    \draw[-,line width=10pt,draw=white] (0,0) .. controls (1,0) and (1,1) .. (2,1);
    \draw (0,0) .. controls (1,0) and (1,1) .. (2,1);
    \draw (0,1) node[left] {$a$};
    \draw (0,0) node[left] {$b$};
    \draw (2,1) node[right] {$b$};
    \draw (2,0) node[right] {$a$};
  \end{tikzpicture}
\end{equation}
The non-equivalence of $\braid{a, b}$ and $\braid{b, a}^{-1}$ corresponds to
the fact that the diagrams (\ref{eq:braid}) and (\ref{eq:braid-inverse}) are not
homotopic relative to their boundaries, as 1-dimensional strings cannot
pass through each other in 3-dimensional space.

If we consider 1-dimensional strings in 4-dimensional space, however, we can
pass the strings through each other by giving one string a higher coordinate
in the 4th dimension so that they do not intersect even when their projections
to other 3 dimensions intersect. This is in fact another invocation of the
Eckmann-Hilton argument, as we ``rotate'' the strings around each other using
the extra dimension. Just as in the 3-dimensional case, there are two ways
to pass these strings through each other up to homotopy, by choosing either $a$
or $b$ to have the higher 4th coordinate.

We can translate this 4-dimensional perspective back to $\omega$-categories
by considering endomorphism 3-cells of
$\apply{\id{\id{\point}}}$. Now there are three dimensions to compose
$a, b : \apply{\id{\id{\point}}} \to \apply{\id{\id{\point}}}$: as
$a \fo_0 b$, as $a \fo_1 b$, or as $a \fo_2 b$. These composites assemble
into a sphere
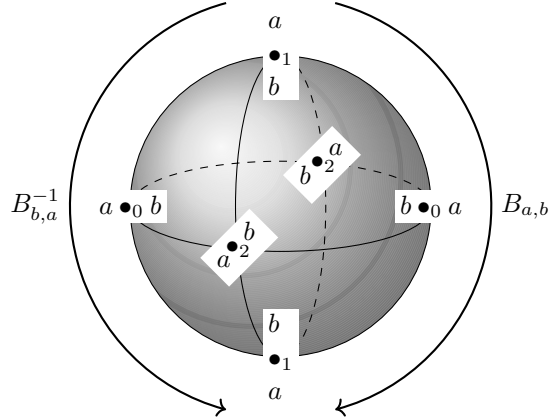
\[
  \begin{tikzpicture}
    \shade[ball color = gray!40, opacity = 0.4] (0,0) circle (2cm);
    \draw (0,0) circle (2cm);
    \draw (-2,0) arc (180:360:2 and 0.6);
    \draw[dashed] (2,0) arc (0:180:2 and 0.6);
    \draw (0,2) arc (90:270:0.6 and 2);
    \draw[dashed] (0,-2) arc (-90:90:0.6 and 2);
    \node [fill=white, inner sep=2pt, align=left] at (-2,0) {$a \fo_0 b$};
    \node [fill=white, inner sep=2pt, align=left] at (2,0) {$b \fo_0 a$};
    \node [fill=white, inner sep=2pt, align=left] at (0,2) {$a$ \\ $\fo_1$ \\ $b$};
    \node [fill=white, inner sep=2pt, align=left] at (0,-2) {$b$ \\ $\fo_1$ \\ $a$};
    \node [fill=white, inner sep=2pt, align=left, rotate=45] at (-0.55,-0.55)
      {\rotatebox[origin=c]{-45}{$a$} \kern -5pt \rotatebox[origin=c]{-45}{$\fo_2$} \kern -5pt \rotatebox[origin=c]{-45}{$b$}};
    \node [fill=white, inner sep=2pt, align=left, rotate=45] at (0.55,0.55)
      {\rotatebox[origin=c]{-45}{$b$} \kern -5pt \rotatebox[origin=c]{-45}{$\fo_2$} \kern -5pt \rotatebox[origin=c]{-45}{$a$}};
    \draw[thick, ->] ([shift={(0,0)}]75:2.8) arc (75:-75:2.8) node[midway, right]{$\braid{a,b}$};
    \draw[thick, ->] ([shift={(0,0)}]105:2.8) arc (105:255:2.8) node[midway, left]{$\braid{b,a}^{-1}$};
  \end{tikzpicture}
\]
i.e.\ there is a 4-cell from any composite to any other composite, including
a generalization of the $\braid{a, b}$ and $\braid{b, a}^{-1}$ cells we
defined above. The two 5-cells $\braid{a, b} \to \braid{b, a}^{-1}$ that
correspond to ``passing strings through each other'' are the front and back
faces of the sphere.
The composites $a \fo_0 b$, $a \fo_1 b$, and $a \fo_2 b$ are all
``the same'' and commutative, as one can be transformed into another by a
rotation of the sphere.

\bibliography{main}{}
\bibliographystyle{plain}

\end{document}